\documentclass{article}
\usepackage[utf8]{inputenc} %%% Included automatically by Overleaf
\usepackage{amsmath}
\usepackage{graphicx} %%% by JHD for inserting figures.

\usepackage{algorithm} %%% by GK 2/9/21 for algorithms
\usepackage{algpseudocode}

\usepackage{xcolor} %%% By GK 2/9/21 for styling text
\usepackage{comment} %%% By JD: allows commenting out large blocks of text
\graphicspath{ {/FIGURES/} }
\usepackage[colorinlistoftodos]{todonotes}%Annotates the text. bordercolor=black, backgroundcolor=orange, textcolor=black linecolor=orange \listoftodos 
\usepackage{dirtytalk}  % Greg 2/10/21
\usepackage{subcaption} % Greg 30/12/21
\usepackage{url}
\usepackage{soul}

\def\ttBH{{\tt BuildHull}}
\def\ttKZCT{{\tt EHD}}
\def\ttHD{{\tt HD}}
\def\ttEHD{{\tt Enhanced Hierarchical Decomposition}}

\def\pps{production possibility set}
\def\PPS{Production Possibility Set}
\def\Pps{Production possibility set}

\def\calA{{\cal A}}
\def\calF{{\cal F}}
\def\calB{{\cal B}}

\def\vrs{\mathop{\tt vrs}\limits}
\def\extB{\mathop{{\tt ext\_}\calB}\limits}

\title{Competing DEA procedures: analysis, testing, and comparisons}
\author{G. Koronakos, J. H. Dul\'a, and D. Despotis}
\date{\today}

\begin{document}
\maketitle

\section*{Abstract}\label{section:Abstract}
Reducing the computational time to process large data sets in Data Envelopment Analysis (DEA) is the objective of many studies. 
Contributions include fundamentally innovative procedures, new or improved preprocessors, and hybridization between -- and among -- all these. 
Ultimately, new contributions are made when the number and size of the LPs solved is somehow reduced.
This paper provides a comprehensive analysis and comparison of two competing procedures to process DEA data sets: \ttBH\ and \ttEHD\ (\ttKZCT).
% We test Khezrimotlagh's recommendation \cite{Khezrimotlagh21}, to use \ttBH\ rather than \ttKZCT\ in sequential (non-parallel) implementations.
A common ground for comparison is made by examining their sequential implementations, applying to both the same preprocessors -- when permitted -- on a suite of data sets widely employed in the computational DEA literature.
In addition to reporting on execution time, we discuss how the data characteristics affect performance and we introduce using the number and size of the LPs solved to better understand performances and explain differences.
Our experiments show that the dominance of \ttBH\ can be substantial in large-scale and high-density datasets. 
Comparing and explaining performance based on the number and size of LPS lays the groundwork for a comparison of the parallel implementations of procedures \ttBH\ and \ttKZCT.

%%%%%%%%%%%%%%%%%%%%%%%%%%%%%%%%%%%%%%%%%%%%%%%%%%%%%%%%%%%%%%%%%%
%%%%%%%%%%%%%%%%%% Begin Section Introduction
%%%%%%%%%%%%%%%%%%%%%%%%%%%%%%%%%%%%%%%%%%%%%%%%%%%%%%%%%%%%%%%%%%
\section{Introduction} \label{section:Introduction}

%%% Start with some  harmless platitudes
%%%
Data Envelopment Analysis (DEA) was introduced in 1978 by Charnes, et al.  \cite{CharnesCooperRhodes78}. 
In DEA, ``Decision Making Units" (DMUs) are engaged in transforming  a common set of $m_1$ inputs into $m_2$ outputs under a specified returns to scale assumption.
Each DMU's measurements for its inputs and outputs constitute a data point with $m=m_1+m_2$ dimensions. 
A constrained linear combination of these points along with ``free-disposability" unit directions of recession defines the \pps\ (PPS).
The type of combination used depends on the returns to scale assumption about the transformation process; here we will focus on  {\it variable returns to scale} (VRS) although the results can be modified and adapted to other DEA returns to scale.
The PPS is a convex, unbounded, polyhedral set.
The data points of the efficient DMUs are on the portion of the boundary of the PPS called the {\it efficient frontier}.
Not all data points on the boundary correspond to efficient DMUs; those not efficient are referred to as {\it weakly efficient}.
All DMUs with data points in the interior of the PPS are inefficient. 
A contribution of the seminal work by Charnes, et al.  \cite{CharnesCooperRhodes78} is a linear programming formulation to calculate the relative efficiency, i.e., the {\it score}, of each of the $n$ DMUs.

\par
A complete DEA study classifies the DMUs as efficient or inefficient.
We refer as the ``standard" approach for the classification and scoring as the one using LPs of the type introduced in \cite{CharnesCooperRhodes78} and \cite{BankerCharnesCooper84}.
These LPs require the solution of $n$ linear programs (LPs), one for each DMU, using all the DMUs' input-output data in the LP's coefficient matrix with at least $m*n$ elements and one or two more rows or columns, depending on the returns to scale assumption. 
Such LPs are referred to as ``full size", which make DEA computationally intensive especially in large scale applications with many DMUs and multiple inputs and outputs.

\par
The computational aspect of DEA is important for accelerating the evaluation process and has been the focus of study since 1980 \cite{BessentBessent80}.
A part of this body of literature is dedicated to preprocessing DEA.
A preprocessor provides information about a DMU's classification or about its score without having to solve an LP.
Preprocesing for DEA is discussed in \cite{SueyoshiChang89},\cite{Sueyoshi90}, \cite{Ali93}, \cite{DulaLopez09},\cite{ChenLai17}, among others.
Preprocessors specifically intended to conclusively identify extreme efficient DMUs are discussed in \cite{Ali93} and \cite{DulaLopez09}.

\par
Many ideas have been proposed for faster DEA algorithms to process increasingly larger DEA data sets.
One of the first approaches based on reducing the size of the LPs to score the DMUs is Ali's ``Reduced Basis Entry" (RBE) \cite{Ali93} procedure.
This procedure excludes the data for inefficient DMUs, as they are identified, in the forward iterations and therefore scores downstream DMUs using smaller LPs.
A full application of RBE will likely finish solving smaller LPs than those at the start especially if the proportion of efficient DMUs, i.e., {\it density}, is low.

\par
Barr \& Durchholz's {\it Hierarchical Decomposition} (\ttHD) \cite{BarrDurchholz97} procedure appeared in 1997 and it applies two phases: in the first phase, data ``blocks" --strict subsets of the data-- are used to identify the inefficient DMUs which are then scored in the second phase using only the efficient DMUs.
Processing the data blocks to identify inefficient DMUs requires smaller LPs and can be done in parallel.
\ttHD\ is perhaps the first DEA procedure that processes VRS hulls of strict subsets of the data to identify inefficient DMUs; we refer to VRS hulls of strict subsets used for this purpose as ``partial" hulls.
A DMU that is inefficient in a partial VRS hull is inefficient for the entire data.
Procedure \ttHD\ sets aside such inefficient DMUs from further consideration in the first phase. 
The survivors are grouped to form the next hierarchy of blocks and the process is repeated until only a single block is derived.
The assessment of this block provides the final block of the efficient DMUs, and marks the end of the first phase of the algorithm.
The final block is used to score all the DMUs in the second phase.
Procedure \ttHD\ requires three user defined parameters, the first to determine the initial block size, the second to increase the block size in the subsequent iterations and the third to terminate the procedure.

\par
Korhonen \& Siitari in \cite{KorhonenSiitari07} apply lexicographic parametric programming to identify the  efficient DMUs in the frontier. 
The algorithm traces a path from unit to unit in the frontier and, when all efficient units are identified, it proceeds to score the other DMUs in a second phase.
This algorithm is enhanced in \cite{KorhonenSiitari09} by halting the process to identify interior points using the hull of a set containing the efficient points.
An important contribution by Chen \& Cho in 2009 \cite{ChenCho09} is a sufficient condition to find the final score of an inefficient DMU prior to access to the full set of efficient DMUs.
This allows scoring of inefficient DMUs in the course of identifying efficient DMUs and potentially accelerates DEA procedures considerably.

\par
Dul\'a proposed the two-phase Procedure \ttBH\ \cite{Dula11} in 2011. 
In the first phase, the extreme points of the \pps\ are identified one at a time using LPs initialized with a few of them -- possibly just one -- until all are found.
The mechanism used to identify new extreme points is hyperplane translation, which requires inner products. 
Therefore, \ttBH\ generates a sequence of nested partial hulls, the elements of which are extreme-points of the full \pps.
The number of LPs solved is deterministic and can be calculated based on how many extreme points are used at the initialization.
The initializing extreme points are discovered using preprocessors.
The size of the LPs at each iteration is given by the number of extreme points of the partial hull at that iteration.
several operations in \ttBH\ can be performed in parallel, such as hyperplane translation which requires inner products. 
Also, as data can be exchanged among the processors, different test points can be processed concurrently by solving LPs in parallel.
Procedures \ttHD\ and \ttBH\ were implemented sequentially and compared in \cite{Dula11}.
In this comparison, \ttBH\ emerged as faster, and in some cases by more than an order of magnitude.

%%%
%%% Gregory, please check this paragraph
%%%
\par
Khezrimotlagh et al. \cite{KhezrimotlaghZhuCookToloo19} introduced in 2019 an enhancement of Procedure \ttHD, since then called \ttEHD\ (\ttKZCT) in Khezrimotlagh \cite{Khezrimotlagh21}.
Procedure \ttKZCT\ also applies two phases with the difference that \ttHD\ is initialized using preprocessors to construct a capacious single block which is used to identify many of the inefficient DMUs in the first phase of the procedure.
Similar to \ttHD, the first phase of \ttKZCT\ partitions the DMUs into two subsets, one containing the DMUs with data points on the boundary of the full VRS \pps\ and the other with the interior (inefficient) DMUs.
Unlike \ttHD, there is not an indeterminate number of levels as new blocks are created from the aggregation of unclassified DMUs and therefore no user specified parameters are required for the cardinality of the blocks at the levels and termination criterion for the first phase in \ttKZCT.
However, \ttKZCT\ requires the user to specify the cardinality of the initializing block.

\par
The boundary points identified at the end of the first phase of \ttKZCT\ may correspond to efficient (extreme and non-extreme) and inefficient (``weak efficient") DMUs. 
In a second phase, these boundary points are used to score the inefficient DMUs.
% Procedure \ttKZCT\ can be directly parallelized; how this is done is discussed in \cite{KhezrimotlaghZhuCookToloo19}.
%\textcolor{red}{In \cite{Khezrimotlagh21} it is claimed that ``\textit{Parallel processing is not applicable for \ttBH\ as its subtasks are dependent on each other}``.}

%% i) \cite{Khezrimotlagh21}: If only one single processor exists to evaluate DMUs, BH could be the fastest available technique, as it requires the least number of calculations in comparison with the existing techniques. However, when the number of processors increases, BH becomes weaker. Parallel processing is not applicable for BH as its subtasks are dependent on each other.

%% ii) Taken from \cite{KhezrimotlaghZhuCookToloo19}: "Clearly, the LPs in BH are not independent from one iteration to another and the algorithm waits until one integration is finished prior to starting the next iteration. In contrast, while BH solves one LP at a time, our procedure allows for solving several LPs at a time, and hence we expect to have a smaller running-time."

\par
In \cite{KhezrimotlaghZhuCookToloo19}, a parallel implementation of \ttKZCT\ was compared with the sequential version of \ttBH\ in terms of execution time for processing DEA data sets. 
It was found that the former was faster than the latter with the right number of parallel processors.
Khezrimotlagh (2021) \cite{Khezrimotlagh21} recommends using \ttBH\ when parallelization is limited or not available, although no testing results are provided.
Jie (2020) \cite{Jie20} proposed a hybridizing \ttBH\ and \ttHD\ where \ttHD\ is parallelized as in \cite{BarrDurchholz97}) and individual blocks are processed using \ttBH.
Jie in \cite{Jie20} reports that his hybridization performs better than his implementation of \ttKZCT\ as parallelized in \cite{KhezrimotlaghZhuCookToloo19}.

\par
%%% Achieving higher performance does not directly imply parallelism, but the primary improvements in highly parallel tasks usually result from enhancing the serial components of the algorithms \cite{Bader08}.
The current work provides a comprehensive comparison of  Procedures \ttKZCT\ and \ttBH\ using a common ground by implementing sequential version of both procedures and  using the same preprocessors and enhancements.
Thus, the current study investigates the speculation made in Section 3 of \cite{Khezrimotlagh21} that \ttBH\ is indicated over \ttKZCT\ in their sequential implementations to process DEA data sets.
We compare the two procedures not only in terms of speed, but also in terms of characteristics that affect their performance; for this purpose, we introduce tracking the number and size of the LPs solved.

\par
The rest of the paper unfolds as follows. 
Section 2 introduces the notation, assumptions and conventions used in this study.
The procedures \ttKZCT\ and \ttBH\ are presented in Section 3, while details about their implementations are given in Section 4.
The procedures are compared in Section 5 by employing the publicly available and widely used {\tt MassiveScaleDEAdata} suite, which consists of 48 highly structured and randomly generated data files.
The procedures are implemented in Matlab and the linear programs (LPs) are solved using Gurobi. 
Finally, conclusions are drawn.

\begin{comment}
Both procedures partition, in a first phase, the DMUs into two sets; one is a subset of the DMUs with data points on the boundary of the \pps, the other subset contains the rest of the DMUs.
In a second phase, the boundary points are used to score the inefficient DMUs.
A difference between the two procedures is the set of boundary points identified at the end of the first phase of \ttKZCT, which can be efficient (extreme and non-extreme) and inefficient (``weak efficient"). 
The set at the end of the first phase of \ttBH\ contains exclusively extreme-efficient DMUs, i.e., the {\it frame}.
\end{comment}

\begin{comment}
%--------FROM 2021 CHAPTER -------------
Khezrimotlagh et al. (2019) improved the HD procedure, called enhanced HD (EHD), and proposed the fastest available technique, in general, to find efficient DMUs as well as inefficiency scores of all DMUs. They only used one computer with several processors to run parallel processing and comparisons; nevertheless, the impact of their approach is noticeable when more than one computer is used in parallel.
%--------------------------
\end{comment}

%%%%%%%%%%%%%%%%%%%%%%%%%%%%%%%%%%%%%%%%%%%%%%%%%%%%%%%%%%%%%%%%%%
%%%%%%%%%%%%%%%%%% Begin Section Notation, Assumptions, and Conventions
%%%%%%%%%%%%%%%%%%%%%%%%%%%%%%%%%%%%%%%%%%%%%%%%%%%%%%%%%%%%%%%%%%
\section{Notation, Assumptions, and Conventions} \label{section:NotationAssumptionsConventions}

Point sets will be denoted using calligraphic letters; e.g. ${\cal A}, {\cal B}$, etc.
Their subsets will be identified using a marking or a superscript: $\bar\calA, {\cal A}^s$, etc.
The cardinality of a set will be expressed as $\vert\cdot\vert$.
For a point set; e.g., ${\cal A}$, with elements in $m$ dimensions and cardinality $n$, we use the uppercase letter, $A$, to denote the $m\times n$ matrix where the points in ${\cal A}$ are the columns of $A$.

\goodbreak 
The following important objects and attributes are used in this study:

\begin{itemize}
\item $\calA=\{a^{1}, \ldots, a^{n}\}$: The full data set of a DEA study.
Each element of $\calA$ consists of the input and output data of a DMU.

\item $n$: the cardinality of $\calA$, $n=\vert\calA\vert$.
\item $m_{1}$: the number of inputs of the DMUs.
\item $m_{2}$: the number of outputs of the DMUs.
\item $m = (m_{1}+m_{2})$: referred to as the \underbar{dimension} of the data points in $\calA$.
The term and parameter applies to a data set referred by its name; i.e., the dimension of $\calA$ is $m$.
\item $\hat m$: number of extreme-efficient DMUs conclusively identified by preprocessors.
\item $X_{i}\in\Re^{m_1}$: the input data vector of DMU $i$, $i= (1, \ldots, n)$.
\item $Y_{i}\in\Re^{m_2}$: the output data vector of DMU $i$, $i= (1, \ldots, n)$.

\item $a^{i}$: the data point for DMU $i$ where:  
\begin{equation}
a^{i}=
\begin{bmatrix}
-X_{i}\\
Y_{i}
\end{bmatrix} \in\Re^m, i= 1, \ldots, n.
\label{eq:TranslatedDataPoint}
\end{equation}

\item $d$: the \underbar{density} of $\calA$. In general, the density of a DEA data set is the proportion of efficient DMUs in the VRS hull of the data. 

\item $\calF$: the \underbar{frame} of $\calA$; the point set containing the extreme points of the VRS \pps\ generated by the DMUs' data; $\calF\subset\calA$. 
The frame $\calF$ is the minimal subset of $\calA$ needed to define the \pps.

\end{itemize}

\noindent The following notation accompanies the presentation and discussion of Procedure \ttKZCT:

\begin{itemize}

\item $\calA^S\subset\calA$: the initializing subset.

\item $p$: the cardinality of the set $\calA^S$: $p=\vert\calA^S\vert$.

\item $\calB^S$: the boundary points of the set $\calA^S$.

\item ${\extB}^S$: the points in $\calA \setminus \calA^S$ which are exterior to $\vrs(\calB^S)$.

\item $\calB$: a point set of boundary points of $\calA$ such that $\calF\subset\calB\subset\calA$. 
The set $\calB$ contains the extreme points of the full \pps\ along with any other efficient or weak efficient DMUs: 

\end{itemize}

\par
As noticed, we work exclusively with the VRS \pps.
The VRS hull defined by a DEA data set, ${\cal A}$, will be referred to as $\vrs({\cal A})$.
For some  ${\cal A}^s \subset {\cal A}$, $\vrs({\cal A}^s)$ will be referred to as a ``partial" hull.
Unless specified, the terms {\it frame}, $\calF$ (without markings), and ``\pps", apply to the full data set, $\calA$ and its VRS hull, and not to any of its subsets or partial hulls.

%%%%%%%%%%%%%%%%%%%%%%%%%%%%%%%%%%%%%%%%%%%%%%%%%%%%%%%%%%%%%%%%%%
%%%%%%%%%%%%%%%%%% Begin Section The two competing DEA procedures
%%%%%%%%%%%%%%%%%%%%%%%%%%%%%%%%%%%%%%%%%%%%%%%%%%%%%%%%%%%%%%%%%%
\section{The competing DEA procedures}\label{section:TheTwoProcedures} \label{section:TheTwoCompetitors}
\par

\subsection{\ttEHD\ (\ttKZCT)} \label{section:ProcedureKZCT}
Procedure \ttKZCT\ is a two-phase procedure proposed by Khezrimotlagh et al. in (2019) \cite{KhezrimotlaghZhuCookToloo19} as an enhancement of \ttHD.
The essential characteristic of the procedure is the use of an initial subset of DMUs in the first phase for the early identification of inefficient DMUs, which are then excluded from the process.
From the remaining unclassified data points, the boundary points of the full \pps\ are obtained.
The inefficient DMUs can then be scored in a second phase.
Procedure \ttKZCT\ is presented next. 

\begin{algorithm}[H]
\renewcommand{\thealgorithm}{}
\floatname{algorithm}{}
\caption{Procedure \ttKZCT\ (Phase 1)}
\textbf{Input:} $\calA$ \Comment{The DEA data set} \\ 
\textbf{Output:} $\calB$  \Comment{Boundary points of $\vrs(\calA)|\calF\subset\calB\subset\calA$}
\begin{algorithmic} [1] 
\State Select a subset $\calA^S$, $\calA^S\subset \calA$.% \newline
\State Identify the points in $\calA^S$ that are on the boundary of $\vrs(\calA^S)$.  
Partition $\calA$ into the boundary points, $\calB^S$, of $\vrs(\calA^S)$ and discard the rest. %\newline
\State Identify the set  ${\extB}^S$: the data points in $\calA \setminus \calA^S$ which are exterior to $\vrs(\calB^S)$.
 %\newline
%Is it better to say "Classify the points to exterior and interior, discard the interior, for the exterior..."
\State Identify the set $\calB$: the data points on the boundary of $\vrs(\calB^S\cup{\extB}^S)$. 
\end{algorithmic}
\end{algorithm}

\par
Procedure \ttKZCT\ uses three different LP formulations is Steps 2, 3, and 4 which remain fixed in size.
Table \ref{table:LPmodelsinKZCT19} summarizes the different LP formulations used in Procedure \ttKZCT. 
Notice that Step 3 employs a ``deleted domain" variant of the standard output oriented VRS LP which permits testing of points not in the LP's coefficient matrix to find the points that are exterior to $\vrs(\calB^S)$.

\begin{table}[H]
\centering
\resizebox{\textwidth}{!}{
\begin{tabular}{lcc}
\hline
\ttKZCT\ & LP model in  Khezrimotlagh, et al. \cite{KhezrimotlaghZhuCookToloo19}\\
\hline
Step 2 & Output oriented VRS Envelopment Form \\
Step 3 &  
$\max \{\varphi: \sum_{j\in \calB^S} x_{j}\lambda_{j} \leq x_{i}, \sum_{j\in \calB^S} y_{j}\lambda_{j} \geq  \varphi y_{jo},  \sum_{j\in \calB^S} \lambda_{j} = 1, \lambda \geq 0\}$
\\
Step 4 & Output oriented VRS Envelopment Form\\ 
\hline
\end{tabular}
}
\caption{LP models employed in Steps 2-4 of \ttKZCT}
\label{table:LPmodelsinKZCT19}
\end{table}

\par
The number and size (in terms of structural variables) of the LPs employed in procedure \ttKZCT\ is given in Table \ref{table:LPcharacteristicsOfKZCT19}.
Both number and size of the LPs solved in this procedure depend on four parameters: the cardinalities of three point sets: $\vert \calA^S\vert = p$, $\vert\calB^S\vert$, and $\vert\calB^S\cup{\extB}^S\vert$ and the fourth being $\hat m$ (the number of extreme elements of the \pps\ identified using preprocessors).

\begin{table}[H]
\centering
%\resizebox{\textwidth}{!}{
\begin{tabular}{lcc}
\hline
\ttKZCT\ &  \begin{tabular}[c]{@{}c@{}}LP size (\# of \\  Structural Variables)\end{tabular} &
\# of LPs \\ 
\hline
Step 2 & $p$  & $p-\hat m$ \\ \\
Step 3 & $|\calB^S|$ & $n-p$ \\ \\
Step 4 & $|\calB^S\cup{\extB}^S|$ & $|\calB^S\cup{\extB}^S|-\hat m$ \\ 
\hline
\end{tabular}
%}
\caption{Characteristics of the optimization process in Steps 2-4 of \ttKZCT}
\label{table:LPcharacteristicsOfKZCT19}
\end{table}

\noindent\underbar{Comments and Observations about Procedure \ttKZCT}

\begin{enumerate}
\item 
Procedure \ttKZCT\ serves as a first phase to process DEA data.
Its output is a set containing the boundary points, or ``best practices", $\calB$, of the full \pps\ of $\calA$.
Since $\calF\subset\calB$, $\calB$ can be used to score the non boundary DMUs in a second phase.

\item  
The cardinality of the initializing subset $\calA^S$ in Step 1 is $p$; it is a user-defined implementation parameter.

\item 
The procedure benefits from a capacious initializing subset $\calA^S$ which will contain many of the \pps's interior points.
This is more likely to happen if $\calA^S$ includes extreme-efficient DMUs and other boundary points.
Preprocessors are available to quickly identify some extreme points of $\vrs(\calA)$ and heuristics have been proposed that will estimate DEA scores potentially identifying additional efficient, or nearly efficient, DMUs.
Such preprocessor and heuristic are advantageous when they require less work than solving linear programs, see \cite{DulaLopez09}.
Let $\hat m$ be the number of extreme-efficient DMUs of the full VRS hull of the data set conclusively identified by preprocessors.
Two preprocessors are employed in \cite{KhezrimotlaghZhuCookToloo19} to derive the subset $\calA^S$ in Step 1 of \ttKZCT.
The first is {\it dimension sorting} \cite{Ali93} that identifies $\hat m$ extreme efficient DMUs such that $1\leq \hat m \leq m$.
The other preprocessor is the pre-scoring heuristic devised in \cite{KhezrimotlaghZhuCookToloo19} to calculate pre-scores for the DMUs depending on their quantile ranking of their input and output values. 
A DMU's pre-score is an estimate of its efficiency score.

\item 
The set $\calB^S\subset \calA^S$ constructed in Step 2 contains the ``best practice" DMUs of $\calA^S$, i.e., the data points of $\calA^S$ on the boundary, of $\vrs(\calA^S)$.
Note that there may be unclassified data points from $\calA \setminus \calA^S$ in the interior and on the boundary of $\vrs(\calA^S)$.
Identifying $\calB^S$ will also distinguish the data points in $\calA^S$ that are interior to the partial hull $\vrs(\calA^S)$, i.e., points that correspond to inefficient DMUs; these are discarded to be scored in a second phase.

%%% and the next set of LPs in Step 4 will also be small.
%%% This situation can occur in low density data sets.

%with $\vert \calA^S \vert -\hat m +1 = p -\hat m +1$ variables (assuming $p\>\hat m+1$) with $m$ constraints can be used to discriminate between those data points on the boundary  and in the interior  of $\vrs(\calA^S)$ needed to create $\calB^S$ in Step 2. The number of such LPs is $\vert \calA^S\vert -\hat m$ where $\hat m$ are the extreme-efficient DMUs identified via a pre-processor.

%% \item \textcolor{red}{ If all extreme-efficient DMUs were to be included in the creation of $\calA^S$ in Step 1, $\calB^S$  will include the full frame $\calF$ at the completion of Step 2, obviating Steps 3 and 4. This is unlikely.}
%%%10/21/21 Removed by JD This ideal situation is highly unlikely given all possible subsets $\calA^S$ of $\cal A$.
%%% The decision on a minimum size of subset $\calA^S that will include the frame is impossible since since the proportion of extreme-efficient DMUs is unknowable.

%As we shall discuss below, the size of the subset also affects the performance of the next steps of the procedure.

\item 
Step 3 classifies the points in $A\setminus \calA^S$ with respect to the partial hull $\vrs(\calB^S)=\vrs(\calA^S)$.
The elements of $A\setminus \calA^S$ can be partitioned into two sets: those that belong to $\vrs(\calB^S)$ not previously identified in Step 2, and those external to $\vrs(\calB^S)$.
The data points in the last category constitute ${\extB}^S$.
The elements of ${\extB}^S$ may be on the boundary or in the interior of the full VRS hull.
Creating ${\extB}^S$ in Step 3 will uncover data points in the interior of the partial hull, $\vrs(\calB^S)$ not in $\calA^S$.
These points can be removed from further consideration until Phase 2.
The total number of the \pps's interior points inside $\vrs(\calB^S)$ is this partial hull's {\it productivity}.
Another category of points that will be identified in Step 3 are data points on the boundary of $\vrs(\calB^S)$ also not in $\calA^S$.
%\textcolor{red}{Strict interpretation of Procedure \ttKZCT\  excludes them from ${\extB}^S$ although, in an implementation this would not affect the effectiveness of the procedure.}
The union of sets ${\extB}^S$ and $\calB^S$ includes the frame $\calF$ of the full hull along with many other points including possibly other of the hull's boundary points.

\item Procedure \ttKZCT\ benefits when $\vert \calB\vert$ is relatively small and the partial hull $\vrs(\calB^S)$ is productive; i.e., it contains a large number of $\calA$'s interior points.
This means that many interior points will be dispatched quickly in Step 3 with small LPs, which also translates to smaller and  fewer LPs in the step that follows.

\item Step 4 processes the elements of $\calB^S\cup{\extB}^S$ to identify the set $\calB$: the ``best practice" DMUs of the full data set $\calA$.
These are data points on the boundary of $\vrs(\calA)$.
The set $\calB$ may include weak efficient or non-extreme efficient DMUs and therefore is a superset of the full frame $\calF$.
The set $\calB$ is used to score the interior DMUs in Phase 2.
%%% An important characteristic of \ttKZCT\ is that the number and size of LPs in Step 4 are almost the same -- the difference being $\hat m$.

%%% \item Moved to section on operation counts
 \item The size  of LPs (number of structural variables) solved in Step 4 is at least the cardinality of the frame.
    This size is augmented by other boundary points in the \pps\ and any of the full hull's interior points not identified in Steps 2 and 3.

\item \label{item:Obviate4thStepOfKZCT}
Step 4 is possible to be obviated.
This occurs when the initializing subset, $\calA^S$,  happens to include all  the boundary points of the full \pps, implying $p\geq\vert\calF\vert$.
If this happens, all DMUs will be scored in Step 3 and there will be neither a need for Step 4 nor for a second phase.
However, if a single boundary point of the \pps\ were to be missing from $\calA^S$, then only the 
$\hat m$ points in $\calB^S \cup {\extB}^S$ are conclusively classified necessitating the fourth step and a second phase.

\item \label{item:RelationsLPsizeAndNumOfLPs}
The information in Table \ref{table:LPcharacteristicsOfKZCT19} permits deductions and insights about  \ttKZCT's workings and performance:

\begin{itemize}
    \item The number of LPs solved in Step 2 and 3 is $n-\hat m$.
     This is the total number of LPs solved by \ttBH, although the sizes differ.
     \item The value $p$ is an upper bound on the size of the LPs in Step 3.
     \item The difference $(p-\vert\calB^S\vert)$ is the number of data points in $\calA^S$ in the interior of its partial hull.
     \item The number of LPs  solved in Step 4, $\vert\calB^S\cup{\extB}^S\vert-\hat m$, and their size   differ by $\hat m$.
     \item The difference between $n$ and the number of LPs solved in Step 4 is the productivity of the partial hull $\vrs(\calA^S)$; i.e., the number of points in the interior of the \pps\ it contains.
      
    \item The size of the frame $\vert\calF\vert$ is a lower bound on the size of the LPs that must be solved in Step 4.
    Also, the size of the frame $\vert\calF\vert$ minus $\hat m$ is a lower bound on the number of the LPs that must be solved in Step 4.
    
\end{itemize}

\end{enumerate}

\subsection{\ttBH}
Procedure \ttBH\ was introduced by J. Dul\'a in 2011 \cite{Dula11}.
It applies to any of the four standard DEA returns to scale, i.e. constant, variable, decreasing and increasing.
For comparison with \ttKZCT\ we confine the discussion to the variable returns to scale.
The key steps of \ttBH\ for the frame identification (first phase) are presented below.

\begin{algorithm}[H]
\textbf{Input:} $\calA$ \Comment{The DEA data set} \\ 
\textbf{Output:} $\calF$  \Comment{The Frame of the DEA data set}
\begin{algorithmic} [1] 
\renewcommand{\thealgorithm}{}
\floatname{algorithm}{}
\caption{Procedure \ttBH\ (Phase 1)}
\State Initialize two dynamic arrays, $\bar\calF ~\&~ \bar\calA$:
%\par \parindent=13pt 
\State $\bar\calF \gets$ any subset of $\calF$. \Comment{$\bar\calF$ is a temporary work-space array}
%\par \parindent=13pt
\State $\bar\calA \gets \calA \setminus \bar\calF$.\Comment{$\bar\calA$ is a temporary work-space array}
\While{$\bar\calA \neq \O$}
% \State 2.2: \begin{minipage}[t]{\linewidth} 
\If{$b\in\vrs(\bar\calF)$}  \Comment{Test point $b$ belongs to the partial hull $\vrs(\bar\calF)$}
        \State $\bar\calA \gets \bar\calA \setminus b$.
        \Else \Comment{Test point $b$ is external to the partial hull $\vrs(\bar\calF)$} 
        \State Find a separating hyperplane between $b$ and $\vrs(\bar\calF)$.
        \State Identify a new frame element, $a^* \in \bar\calA$, using the separating hyperplane.
             \par       %        \parindent=33pt %
         \State $\bar\calF \gets \bar\calF \cup a^*$;
          \State $\bar\calA \gets \bar\calA \setminus a^*$.
        \EndIf
\EndWhile  
\State Conclusion: $\calF \gets \bar\calF$.

\end{algorithmic}
\end{algorithm}

\par
The LP model formulated by Dul\'a (2011) \cite{Dula11} to test whether a point belongs to the partial hull of $\vrs(\bar\calF)$ in Step 5 of \ttBH\ is given below: 

\[
 \left \{\min_{\delta \geq 0, \lambda \geq 0} \delta: \sum_{i\in \bar\calF} a^{i}\lambda_{i} + e\delta \geq b, \sum_{i\in \bar\calF} \lambda_{i} = 1 \right\} \tag{3} 
 \label{eq:ModelBH}
\]
\par
The dual of program (\ref{eq:ModelBH}) provides the separating hyperplane for an external test point.
The auxiliary variable, $\delta$, estimates the separation from the VRS hull of the points in $\bar \calF$ to the test point $b$.
The vectors $a^i;\; i=1,\ldots,n$ are as defined in Expression (\ref{eq:TranslatedDataPoint}) and $e$ is a vector of 1s.

\bigskip
\goodbreak
\noindent\underbar{Comments and Observations about Procedure \ttBH}

\begin{enumerate}
\item
Procedure \ttBH\ generates a sequence of nested subsets of $\calF$ which grows one frame element at a time, terminating with the full frame.

\item In \cite{Dula11} Procedure \ttBH\ is illustrated by employing a single extreme point of the \pps\ for the initialization in Step 2. 
However, as noted by Dul\'a (2011) \cite{Dula11}, any number of initializing extreme elements of the hull can be used in this step; e.g., as could be obtained with dimension sorting.

%%% Almost verbatim from my paper.
%%% 10/3/21 Excerpt from IJoC2011:
%%% "Any known extreme element of a hull can be inserted in the initializing partial hull."

\item 
Each iteration in the {\tt while} loop in Step 4 will reduce the cardinality of the set $\bar\calA$ by one. 
This means that Procedure \ttBH\ will perform exactly $n-|\bar\calF|$ iterations, where $|\bar\calF|$ is the cardinality of the initializing partial frame.
%%% This makes Procedure \ttBH\ fully deterministic.

\item
Step 5 begins by determining whether a test point belongs to the current partial hull or is external.
This can be carried out by solving an LP of dimension $m$ by the number of extreme points in the current partial hull.
Such an LP must also provide a separating hyperplane between the partial hull, $\vrs(\bar\calA)$ and the test point, $b$.
For this purpose, model (\ref{eq:ModelBH}), originally proposed in \cite{Dula11}, can be employed.

\item
A test point, $b$, that belongs to a partial hull may be in its interior or on the boundary as an non-extreme efficient DMU or weak efficient. 
Therefore, Procedure \ttBH\ permits only extreme-efficient DMUs to belong to $\bar\calF$.
Thus, in contrast to the output of first phase of \ttKZCT\ that may contains ``weak efficient" DMUs, $\calF$ will include only extreme-points of the full VRS hull upon completion, i.e., the {\it frame}.

\item
If the test point, $b$, is exterior to the partial hull $\bar\calF$ in Step 7, then there exists a hyperplane that separates it from the partial hull.
Translating this separating hyperplane away from the partial hull will necessarily uncover a new frame element, $a^*\in\bar\calA$.
This operation involves sorting simple inner products between the hyperplane's parameters and the points in $\bar\calA$.
A newly uncovered extreme point, $a^*$, in Step 9 is not necessarily the same as the test point $b$.
The possibility of ties when translating the separating hyperplane -- and how to resolve them to locate a single new extreme point -- is discussed in \cite{Dula11}.

\item 
The number of iterations which require hyperplane translation is $\vert\calF\vert-\hat m$, while the number of hyperplane translations, i.e., number of inner products, performed in each one of these iterations is specified by the cardinality of the unclassified points at that iteration.
\end{enumerate}

%In each next iteration the number of hyperplane translations will be decreased by one.

\par
For Procedure \ttBH, the number of LPs solved is deterministic at $n-\hat m$, where $\hat m$ is the initializing partial frame at the beginning of the procedure.
On the other hand, in \ttKZCT\ only the Steps 2-3 require $n-\hat m$ LPs (Table \ref{table:LPcharacteristicsOfKZCT19}).
The size of the LPs solved in \ttBH\ is determined by the number of frame elements identified at a given iteration as the procedure progresses. 
The LPs remain fixed at $\vert \calF \vert$ structural variables after the iteration in which all the frame elements are identified.

\subsection{Impact on Phase 2} \label{section:ImpactPhaseTwo}
Inefficient DMUs are scored using the output of the procedures, $\calB$ for \ttKZCT\ and $\calF$ for \ttBH, in a second phase.
The DMUs in these two sets will define the structural variables of the DEA LPs used in the scoring phase.
A difference between \ttKZCT\ and \ttBH\ is that $|\calB|\geq|\calF|$ since the former includes efficient non-extreme and weak efficient DMUs.
This difference can affect the performance in Phase 2.  %%% This is discussed elsewhere.
Phase 2 may have to solve larger LPs using $\calB$ after \ttKZCT\  but the number of DMUs to be scored will be fewer and vice-versa after using \ttBH.

%%%%%%%%%%%%%%%%%%%%%%%%%%%%%%%%%%%%%%%%%%%%%%%%%%%%%%%%%%%%%%%%%%%%%%%
%%%%%%%%%%%%%%%%%  Begin section TheTwoImplementations Section %%%%%%%%%%%%%%%%
%%%%%%%%%%%%%%%%%%%%%%%%%%%%%%%%%%%%%%%%%%%%%%%%%%%%%%%%%%%%%%%%%%%%%%%
\section{Implementations of \ttKZCT\ \& \ttBH\ }\label{section:TheTwoImplementations}
The implementations of \ttKZCT\ and \ttBH\ along with the characteristics of the computational environment used in \cite{KhezrimotlaghZhuCookToloo19} and in the current study, are presented in this section.

\subsection{Implementation of \ttKZCT\ \& \ttBH\ in \cite{KhezrimotlaghZhuCookToloo19}} \label{section:OurImplementationsOfBothProcedures}
A fundamental feature of Procedure \ttKZCT\ is the cardinality, $p$, of subset $\calA^S$.
The implementation of \ttKZCT\ in \cite{KhezrimotlaghZhuCookToloo19} sets this parameter to $p=\sqrt{n}$ before any step is taken. 
Two methods are employed in \cite{KhezrimotlaghZhuCookToloo19} to select elements for $\calA^S$ in Step 1  \ttKZCT.
The first is the ``dimension sorting" preprocessor, discussed above, which identifies $1\leq\hat m\leq m$ extreme points of the \pps.
The second method is a heuristic for calculating a ``pre-score" for each DMU to serve as a proxy of its efficiency score, which was introduced in \cite{KhezrimotlaghZhuCookToloo19}. 
Although the method cannot conclusively identify efficient DMUs, the effort involved is relatively low.
The objective of the two preprocessors used in \cite{KhezrimotlaghZhuCookToloo19} when constructing $\calA^S$ is to include data points of DMUs known to be on, or near, the efficient frontier.

\par
Procedure \ttBH\ is implemented fully sequentially (not parallelized) in \cite{KhezrimotlaghZhuCookToloo19}. 
In particular, in \cite{KhezrimotlaghZhuCookToloo19}, Procedure \ttBH\ is initialized in Step 1  with a single element of the frame, $\calF$.
Hence, based on this implementation of \ttBH, a total of $n-1$ LPs are solved to identify the full frame of $\vrs(\calA)$.

\subsection{Implementation of \ttKZCT\ \& \ttBH\ in the current study} \label{OurImplementationsOfBothProcedures}

The implementation of Procedure \ttKZCT\ in this study follows the fully sequential (non-parallelized) version in \cite{KhezrimotlaghZhuCookToloo19} as presented above.  
The implementation for Procedure \ttBH\ differs from \cite{KhezrimotlaghZhuCookToloo19} in that it incorporates the two preprocessing methods used also for \ttKZCT\, i.e., i) Ali's \cite{Ali93}``dimension sorting" to find $\hat m\leq m$ extreme points of the \pps\ and ii) the pre-scoring heuristic.

\par
Dimension sorting applied to \ttBH\ reduces the number of LPs solved from $n-1$ to $n-\hat m$. 
The pre-scores obtained using the pre-scoring heuristic provide an estimate of the DMUs' proximity to the efficient frontier.
These pre-scores are used in our implementation of \ttBH\ to reorder $\calA$ in ascending order upon entering Step 4 of the procedure.
The effect of this reordering of the data points is anticipated to reduce the average LP size by attenuating the rate at which the LPs grow.

\par
Table \ref{table:BHImplementations} presents the characteristics of the implementations of \ttBH\ in \cite{KhezrimotlaghZhuCookToloo19} and in the current study.

\begin{table}[H]
\centering
\resizebox{\textwidth}{!}{
\begin{tabular}{lcccc}
\hline
\begin{tabular}[c]{@{}c@{}}Implementation \\ of \ttBH\end{tabular} &
  LP Formulation &
  \begin{tabular}[c]{@{}c@{}}Size of initial \\ partial hull $\bar\calF$ \end{tabular} &
  \begin{tabular}[c]{@{}c@{}} Total\\ \# of LPs \\ \end{tabular} &
  \multicolumn{1}{c}{\begin{tabular}[c]{@{}c@{}}Order of processing \\ DMUs \end{tabular}} \\ \hline
\begin{tabular}[c]{@{}c@{}}Khezrimotlagh \\ et al. (2019) \cite{KhezrimotlaghZhuCookToloo19}\end{tabular} &
  Model (\ref{eq:ModelBH}) & $1$ & $n-1$ & Random \\
Current study &  Model (\ref{eq:ModelBH}) & $ \hat m$ & $n-\hat m$ &
  \begin{tabular}[c]{@{}c@{}} Ascending order\\ of their prescores\end{tabular} \\ 
  \hline
\end{tabular}
}
\caption{Implementation of \ttBH\ in \cite{KhezrimotlaghZhuCookToloo19} compared to the current study}
    \label{table:BHImplementations}
\end{table}

\par
Based on Tables \ref{table:LPcharacteristicsOfKZCT19} and \ref{table:BHImplementations} a comparison of the number of LPs solved in \ttKZCT\ and \ttBH\ can be carried out.
As can be seen, the number of LPs solved in Steps 2-3 of \ttKZCT\ is the same as the  number of LPs solved in \ttBH\ in total.
Procedure \ttKZCT\ solves additional LPs in Step 4.

\subsection{Implementation Environment} \label{Enviroment}

The environment used for the two implementations in this study resembles the corresponding one used in \cite{KhezrimotlaghZhuCookToloo19}.
We differ in the optimization solver employed: we use Gurobi\footnote{gurobi.com/products/gurobi-optimizer.} rather than Matlab's {\tt linprog()}\footnote{mathworks.com/help/optim/ug/linprog.html.}.
Gurobi offers both the primal and dual simplex algorithms to solve LPs; both are used in this study.
The default LP algorithm used in Matlab's {\tt linprog()} is the dual simplex; the primal simplex is not an option.
The  specifications of the computational environment used in \cite{KhezrimotlaghZhuCookToloo19} and the current study are given in Table \ref{table:HardwareSpecs} below.

\begin{table}[H]
  %\label{table:HardwareSpecs} %%% Check this label.
\centering
\resizebox{\textwidth}{!}{
\begin{tabular}{lcc}
\hline
Environment & Khezrimotlagh et al. \cite{KhezrimotlaghZhuCookToloo19} & Current Study \\ 
\hline
Operating System &  Windows 10 64-bit &  Windows 10 64-bit \\
Processor & Intel®Core™ i7-7820HK @2.90GHz  & Intel®Core™ i7-8665U @1.9GHz \\ 
Random Access Memory (RAM) &  16GB &  16GB \\
Scripting Language & Matlab 2017a & Matlab R2019b \\
Optimization Solver & Matlab linprog & Gurobi v9.1.2  \\
\hline
\end{tabular}
}
\caption{Computational environment specifications in \cite{KhezrimotlaghZhuCookToloo19} and in the current study}
\label{table:HardwareSpecs}
\end{table}

%%%%%%%%%%%%%%%%%%%%%%%%%%%%%%%. Begin Section "Comparison and Results" %%%%%%%%%%%%%%%%%%%%%%%%%
\section{Results and Comparison of Procedures \ttKZCT\ \& \ttBH}
\label{section:ComparisonOfTwoProceduresWithOurImplementation}

The comparison between \ttKZCT\ and \ttBH\ is made using a publicly available\footnote{\url{dsslab.cs.unipi.gr/datasets/MassiveScaleDEAdata}.} {\tt MassiveScaleDEAdata} suite of 48 synthetic data files created by J. Dul\'a.
The datasets has been already widely employed in the literature, e.g., in \cite{KorhonenSiitari07}, \cite{Dula08}, \cite{KorhonenSiitari09}, \cite{DulaLopez09}, \cite{Dula11} and \cite{Jie20}, among others.
Detailed information about the data generation process is provided in \cite{Dula08}.

\par
The two Procedures \ttKZCT\ and \ttBH\ are compared with respect to the three characteristics of a DEA data set -- cardinality $n$, dimension $m$, and density $d$.
The data sets vary in  four  cardinality  levels  ($n =$ 25K, 50K, 75K, 100K), four dimensions ($m = 5, 10, 15, 20$) and three density levels ($d = 1\%, 10\%, 25\%$).
For six specific data sets the number of frame elements differs slightly from the nominal density: {\tt 10by25000at10} has $|\calF|=2524$; {\tt 15by25000at01} has $|\calF|=262$ {\tt 20by25000at01} has $|\calF|=252$; {\tt 10by50000at10} has $|\calF|=5004$, {\tt 15by50000at01} with $|\calF|=501$; and {\tt 10by75000at10} has $|\calF|=7501$.

\par 
None of the synthetic DEA data sets of the {\tt MassiveScaleDEAdata} suite utilized in this study contain weak efficient or non-extreme efficient DMUs.
The boundary points of the VRS \pps\ are limited to the data of extreme-efficient DMUs.
This means both procedures produce the same output,  $\calB=\calF$, in the 48 instances.
Therefore, there is no need to consider a second phase in the comparison.
Finally, in all 48 instances of the data suite, dimension sorting identified the maximum  number of extreme elements of the \pps, i.e. $\hat m= m$.

\par
The summary statistics about the execution time for the 48 data sets employing Gurobi's primal and dual simplex algorithms in both procedures are presented in Table \ref{table:TimePrimalDual}.
Based on the average times to process the 48 data sets, \ttKZCT\ benefits more from the Gurobi's primal simplex than \ttBH.
Therefore, we report on the experimental results of our implementations of \ttKZCT\ and \ttBH\ using Gurobi's primal simplex algorithm.

\begin{comment}

\begin{table}[H]
  \centering
  \resizebox{\textwidth}{!}{
    \begin{tabular}{lcccc}
    \hline
    \begin{tabular}[c]{@{}c@{}}Time \\(in seconds) \end{tabular}&
    \begin{tabular}[c]{@{}c@{}}\ttKZCT\ \\Dual Simplex \end{tabular} &
    \begin{tabular}[c]{@{}c@{}}\ttKZCT\ \\Primal Simplex \end{tabular} &
    \begin{tabular}[c]{@{}c@{}}\ttBH\ \\Dual Simplex \end{tabular} &
    \begin{tabular}[c]{@{}c@{}}\ttBH\ \\Primal Simplex \end{tabular} \\ 
    \hline
   \begin{tabular}[c]{@{}}Min  \end{tabular} & 39.86 & 40.70 & 31.67 & 39.18 \\[0.06cm]
   \begin{tabular}[c]{@{}}Max  \end{tabular}  & 5372.21 & 5051.58 & 2123.27 & 2295.94 \\[0.06cm]
   \begin{tabular}[c]{@{}}Mean \end{tabular} & 936.62 & 851.66 & 368.72 & 415.39 \\
   \begin{tabular}[c]{@{}}Standard \\Deviation \end{tabular}& 1191.24 & 1223.88 & 466.64 & 509.78 \\
   \begin{tabular}[c]{@{}}Coefficient \\of Variation\end{tabular} & 1.27 & 1.44 & 1.27 & 1.23 \\
   \hline
   \end{tabular} 
}
    \caption{Descriptive statistics of execution time for the 48 data sets using Gurobi's Primal and Dual Simplex}
  \label{table:TimePrimalDual}  
\end{table}

\end{comment}

\begin{table}[H]
  \centering
  \resizebox{\textwidth}{!}{
    \begin{tabular}{lcccc}
    \hline
    \begin{tabular}[c]{@{}c@{}}Time \\(in seconds) \end{tabular} &
    \begin{tabular}[c]{@{}c@{}}\ttKZCT\ \\Dual Simplex \end{tabular} &
    \begin{tabular}[c]{@{}c@{}}\ttKZCT\ \\Primal Simplex \end{tabular} &
    \begin{tabular}[c]{@{}c@{}}\ttBH\ \\Dual Simplex \end{tabular} &
    \begin{tabular}[c]{@{}c@{}}\ttBH\ \\Primal Simplex \end{tabular} \\
    \hline
    \begin{tabular}[c]{@{}c@{}}Min \end{tabular} & 39.86 & 40.70 & 31.67 & 39.18 \\[0.06cm]
    \begin{tabular}[c]{@{}c@{}}Max \end{tabular} & 5372.21 & 5051.58 & 2123.27 & 2295.94 \\[0.06cm]
    \begin{tabular}[c]{@{}c@{}}Mean \end{tabular} & 936.62 & 851.66 & 368.72 & 415.39 \\
    \begin{tabular}[c]{@{}c@{}}Standard \\Deviation \end{tabular} & 1191.24 & 1223.88 & 466.64 & 509.78 \\
    \begin{tabular}[c]{@{}c@{}}Coefficient \\of Variation \end{tabular} & 1.27 & 1.44 & 1.27 & 1.23 \\
    \hline
    \end{tabular}
  }
  \caption{Descriptive statistics of execution time for the 48 data sets using Gurobi's Primal and Dual Simplex}
  \label{table:TimePrimalDual}
\end{table}

\par
\begin{figure}[H]
\begin{center}
  \includegraphics[width=1.00\linewidth,height=0.55\linewidth]{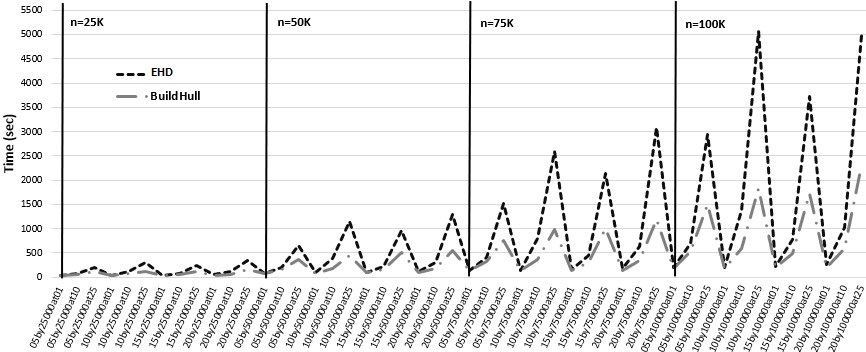}
  \caption{Execution times for all data sets in the {\tt MassiveScaleDEAdata} sorted by Cardinality$\rightarrow$Dimension$\rightarrow$Density.}
  \label{fig:RunningTimes}
\end{center}
\end{figure}

The experimental results also show that the average time for hyperplane translation for \ttBH\
in the 48 data sets of {\tt MassiveScaleDEAdata} suite is less than 1.5\% of the total running time.
This is hardly enough to impact the comparisons, thus is omitted from the analyses.
Figure \ref{fig:RunningTimes} summarizes the performance of the two procedures in terms of running times for the 48 data sets in the {\tt MassiveScaleDEAdata} suite.
The detailed results, organized by cardinality, are given in the appendix, in Tables \ref{table:KZCT19Steps25K}-\ref{table:KZCT19Steps100K} for \ttKZCT\ and in Tables \ref{table:BestWorstBH25K}-\ref{table:BestWorstBH100K} for \ttBH.
These tables include the total running time along with the number and size of the solved LPs, stage-wise for \ttKZCT\ and directly for \ttBH.

\par
The results allow the following conclusions about the two procedures' running times as evident in  Figure \ref{fig:RunningTimes}:

\begin{enumerate}
    \item \ttBH\ is faster than \ttKZCT\ in all 48 instances for the data sets in the {\tt MassiveScaleDEAdata} suite using sequential implementations.
    The advantage of \ttBH\ over \ttKZCT\ can be attributed to the additional large LPs solved in Step 4 of the latter.
    \item The procedures are closest in time for {\tt 05by25000at01} with 39.74 secs.~for \ttBH\ vs 40.70 secs.~for \ttKZCT, while the biggest difference occurs at {\tt 10by100000at25} with 1,817.77 secs.~for \ttBH\ vs 5,051.58 secs.~for \ttKZCT.
    \item For a given dimension and density, running time increases with cardinality.
   \item Within each cardinality we can appreciate how time clearly increases with density and how the effect is more severe for Procedure \ttKZCT.
    \item The closest gaps between the two procedures occur at $d=1\%$.
\end{enumerate}

\par 
In the following discussion, computational results are presented and discussed  in terms of the effect of density, $d$, cardinality, $n$, and dimension, $m$, using representative values of the parameters.

\begin{figure}[H]
\begin{subfigure}{.5\textwidth}
  \centering
  \includegraphics[width=1\linewidth, height=0.75\textwidth]{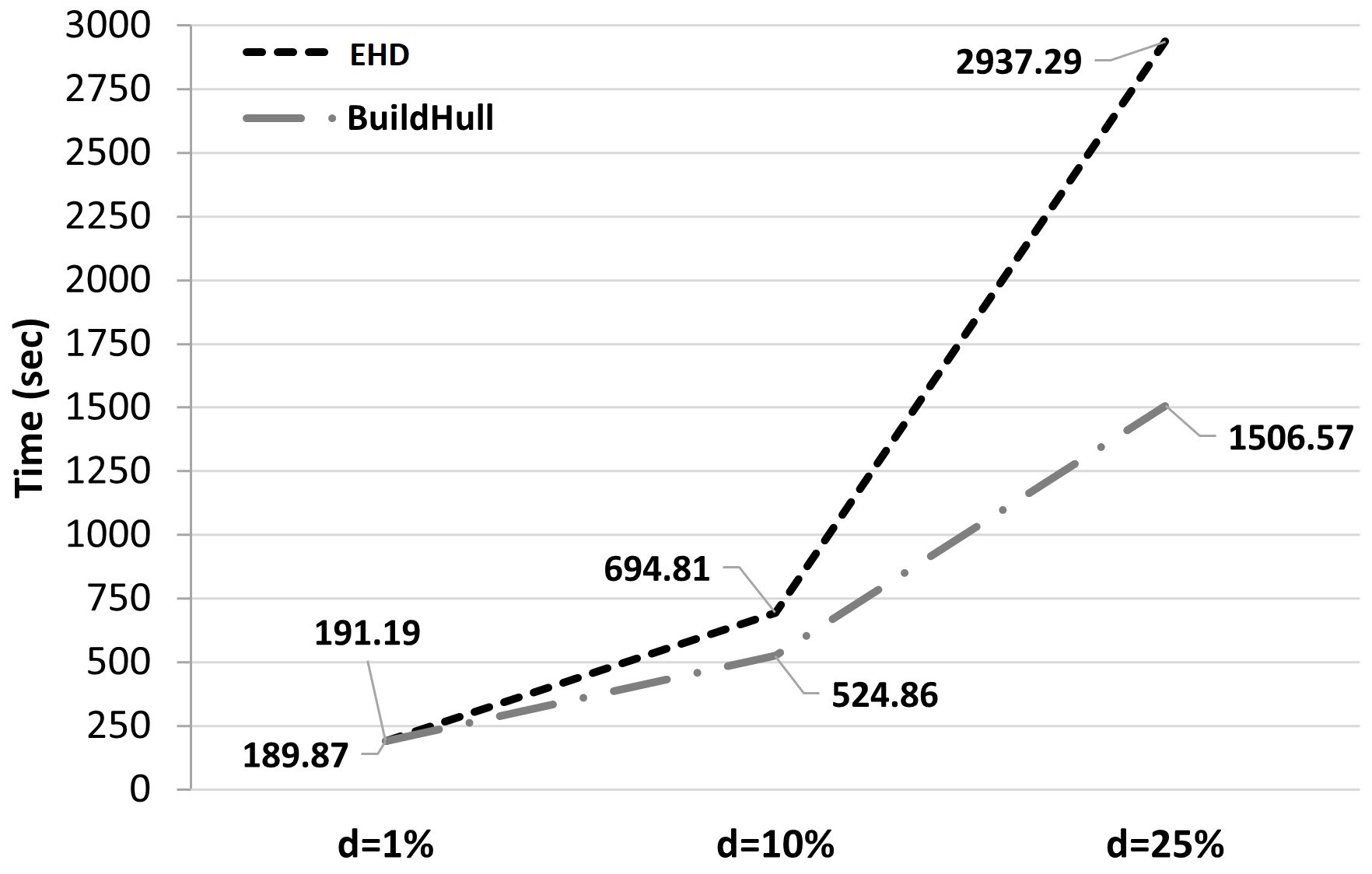}  
  \caption{$n=100$K, $m=5$}
  \label{fig:fig2a}
\end{subfigure}
\begin{subfigure}{.5\textwidth}
  \centering
  \includegraphics[width=1\linewidth, height=0.75\textwidth]{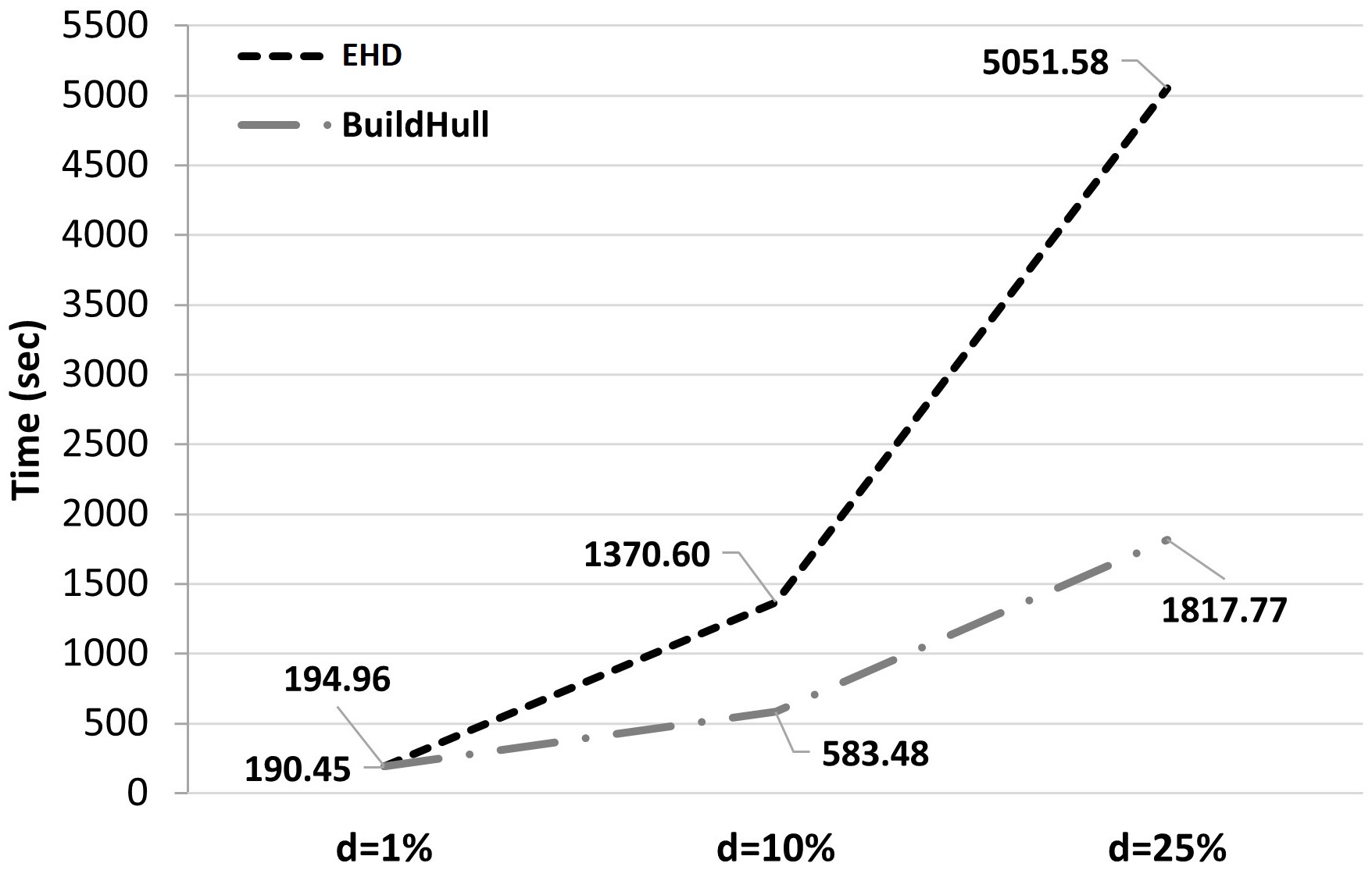}  
  \caption{$n=100$K, $m=10$}
  \label{fig:fig2b}
\end{subfigure}
\begin{subfigure}{.5\textwidth}
  \centering
  \includegraphics[width=1\linewidth, height=0.75\textwidth]{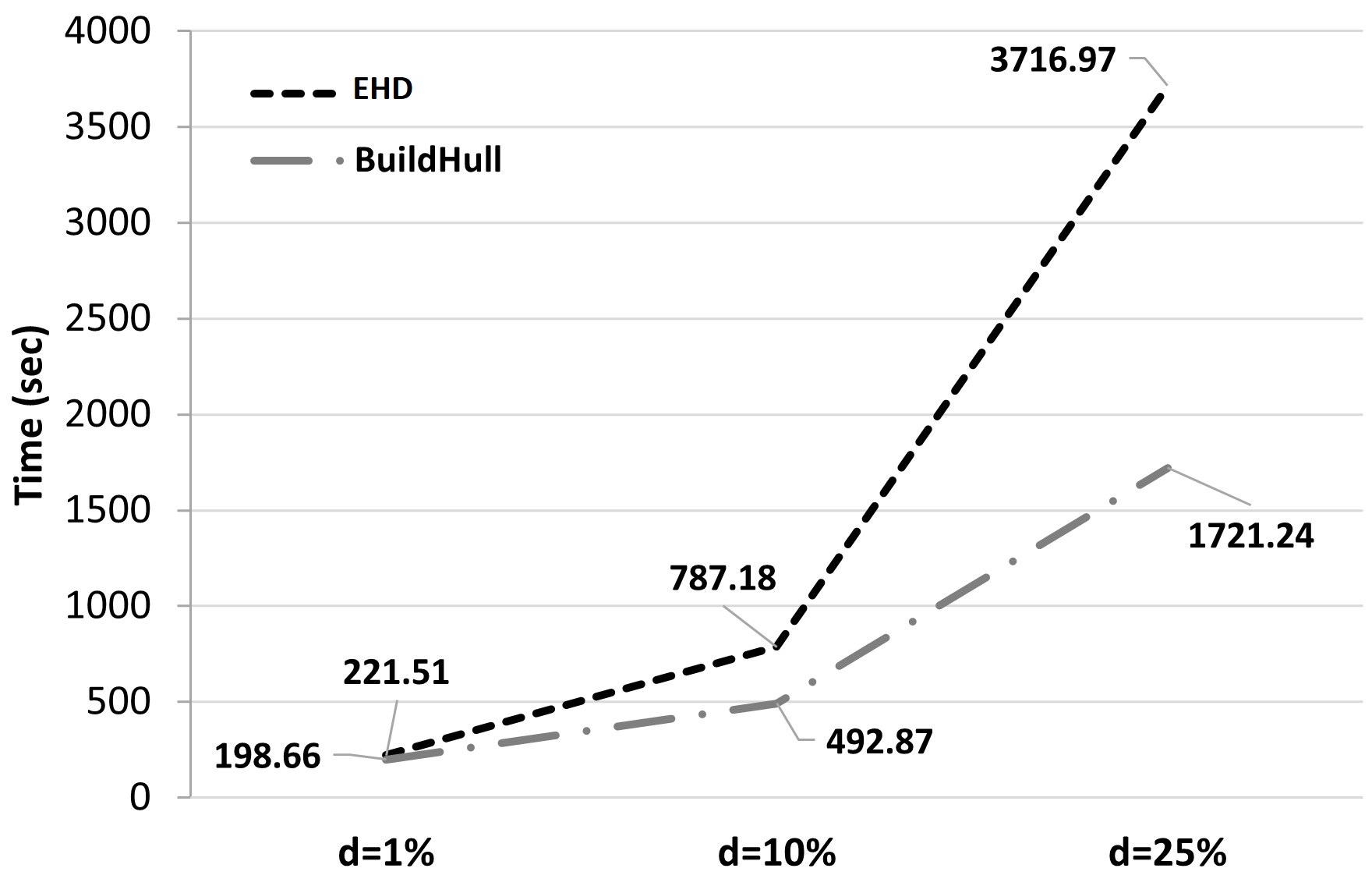}  
  \caption{$n=100$K, $m=15$}
  \label{fig:fig2c}
\end{subfigure}
\begin{subfigure}{.5\textwidth}
  \centering
  \includegraphics[width=1\linewidth, height=0.75\textwidth]{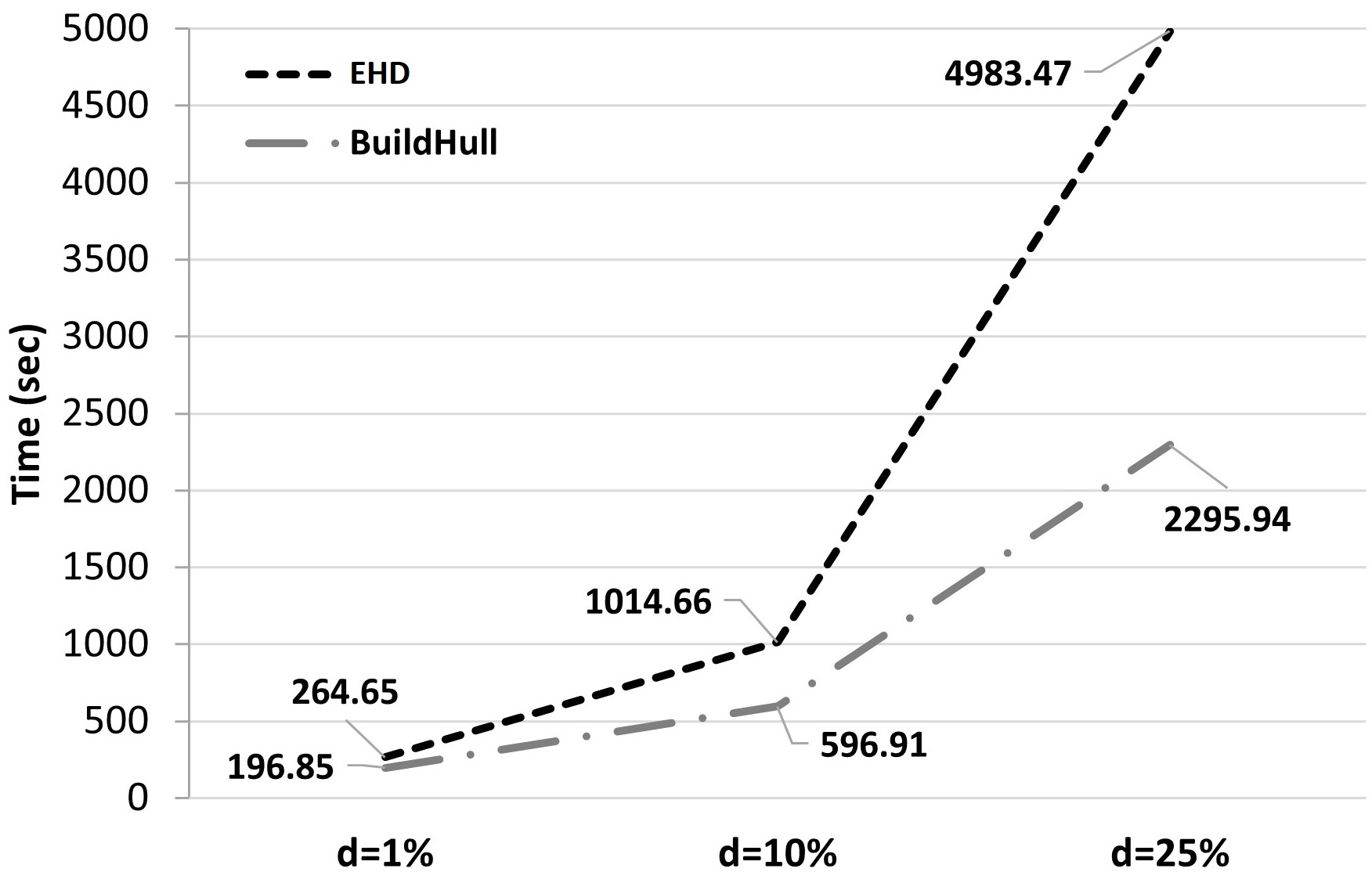}  
  \caption{$n=100$K, $m=20$}
  \label{fig:fig2d}
\end{subfigure}
\caption{Impact of density on running time for the data sets in the {\tt MassiveScaleDEAdata} suite when the cardinality is $n=$100K.}
\label{fig:Density100k}
\end{figure}

\bigskip
\goodbreak
%\textcolor{red}{Impact of density}
\noindent \underbar{Density: $d$}. The four panels in Figure \ref{fig:Density100k} compare the effect of density, $d$, on the running time performance of the two procedures  for the data sets with highest cardinality, $n=$100K, across four dimensions.
The four plots depict how the two procedures start at almost the same speed when $d=1\%$ and, as the density increases,  \ttBH\ gains an advantage over \ttKZCT\ ending up  more than twice as fast  when $d=25$.
This relative performance based on density using these massive data sets is representative of what occurs with the rest of the point sets in the {\tt MassiveScaleDEAdata} suite.

\par
As noticed, the number of LPs solved in Steps 2 and 3 of \ttKZCT\ equals the total number of LPs solved in \ttBH. 
Also, Step 2 solves $(p-\hat m)$ LPs of size  $p=\sqrt{n}$, while Step 3 solves $(n-p)$ LPs of size less than or equal to $p$, see Table \ref{table:LPcharacteristicsOfKZCT19}.
Given that there are only frame elements on the boundary of the \pps\ in the {\tt MassiveScaleDEAdata} suite and that $\hat m=m$ in all its instances, a lower bound on the average LP size in Step 4 is $\vert \calF\vert$.
When both density and cardinality are small, $\vert \calF\vert$ is not much greater than $p$: e.g. $d=1\%, n=25K \Rightarrow \vert \calF\vert=250, p=\sqrt{25,000}=158$.
However, as $d$ and $n$ increase, the difference between  $\vert \calF\vert$ and $p$ becomes much greater to the point when $d=25\%, n=100K \Rightarrow \vert \calF\vert=25,000, p=316$.
This explains the relatively narrow running time differences for the two procedures when density and cardinality are small and the large gap when $\vert \calF\vert \gg \sqrt{n}$ which occurs when density is high in large cardinalities portrayed in Figures \ref{fig:RunningTimes} and \ref{fig:Density100k}.

\begin{figure}[H]
\begin{subfigure}{.5\textwidth}
  \centering
  \includegraphics[width=1\linewidth, height=0.75\textwidth]{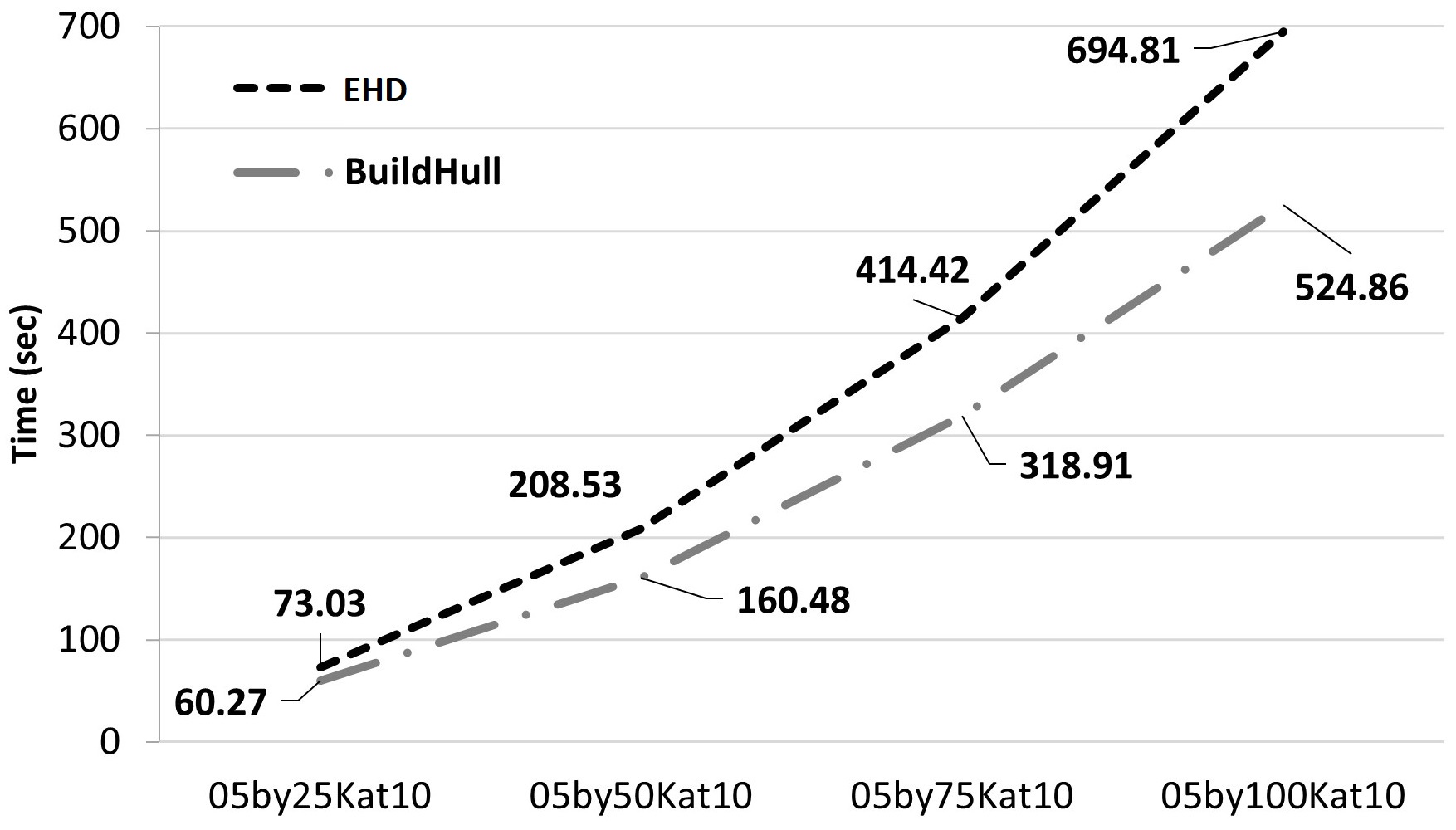}  
  \caption{$d=10\%$, $m=5$}
  %\label{fig:sub-first}
\end{subfigure}
\begin{subfigure}{.5\textwidth}
  \centering
  \includegraphics[width=1\linewidth, height=0.75\textwidth]{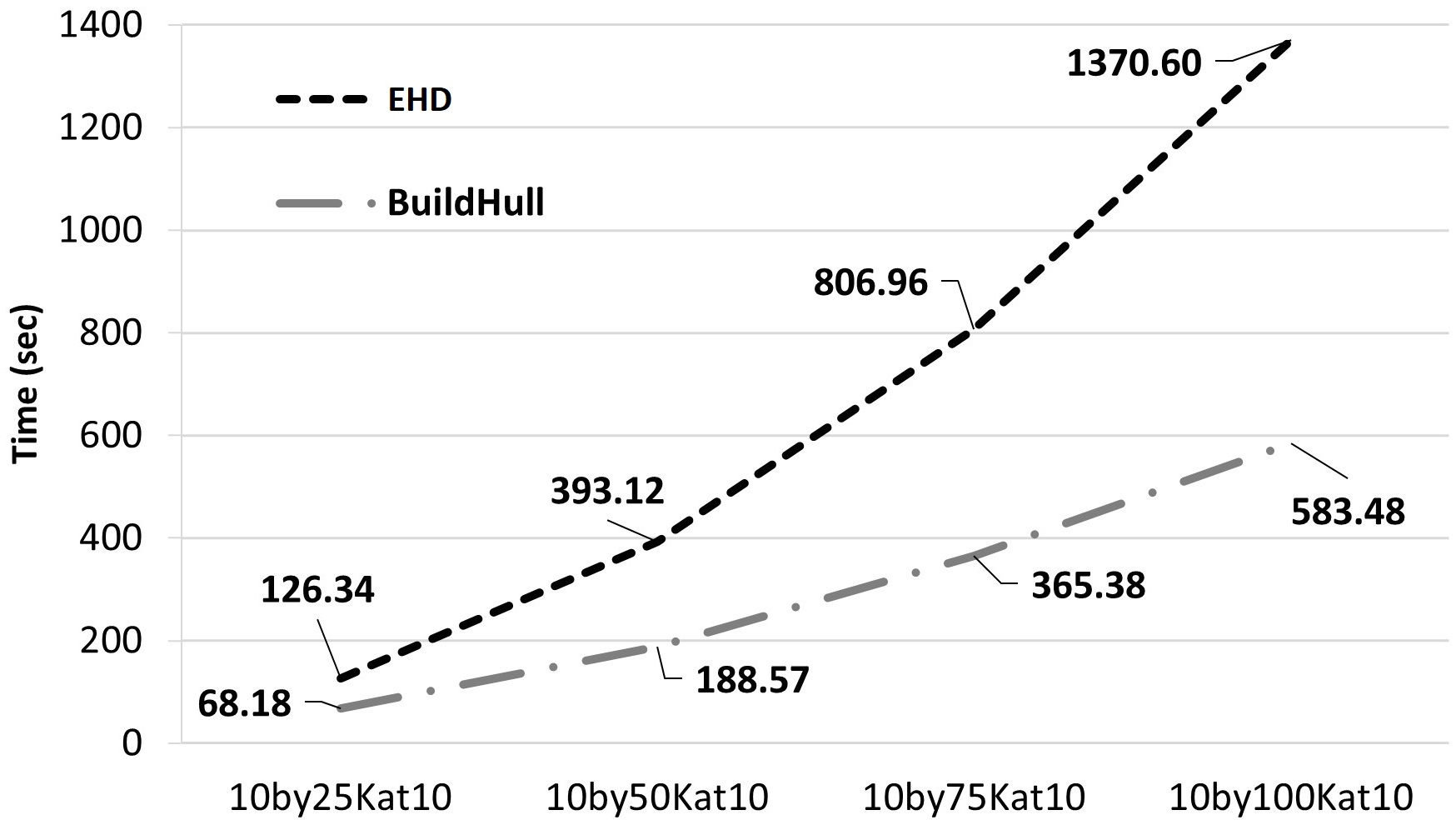}  
  \caption{$d=10\%$, $m=10$}
  %\label{fig:sub-second}
\end{subfigure}
\begin{subfigure}{.5\textwidth}
  \centering
  \includegraphics[width=1\linewidth, height=0.75\textwidth]{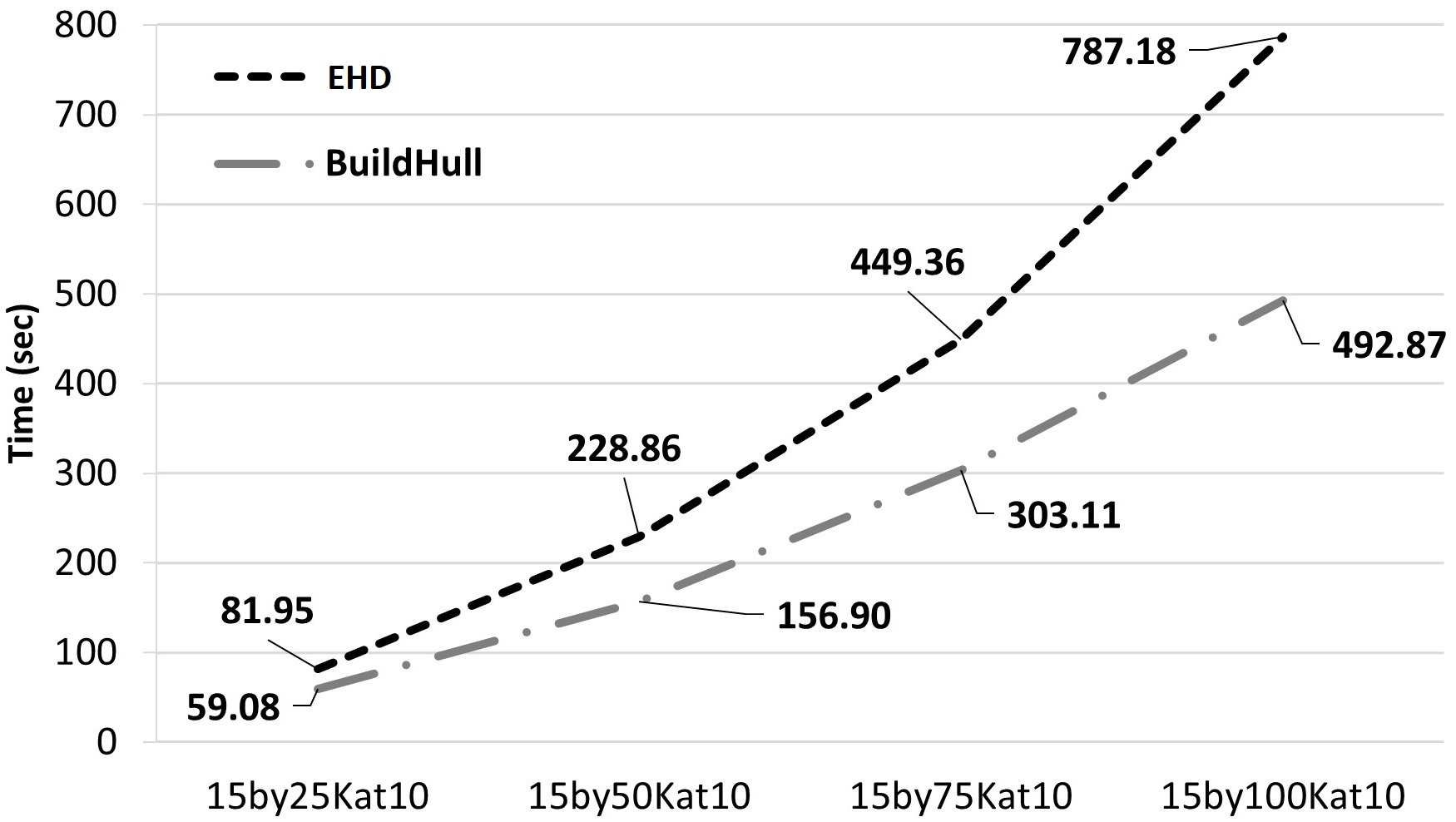}  
  \caption{$d=10\%$, $m=15$}
  %\label{fig:sub-second}
\end{subfigure}
\begin{subfigure}{.5\textwidth}
  \centering
  \includegraphics[width=1\linewidth, height=0.75\textwidth]{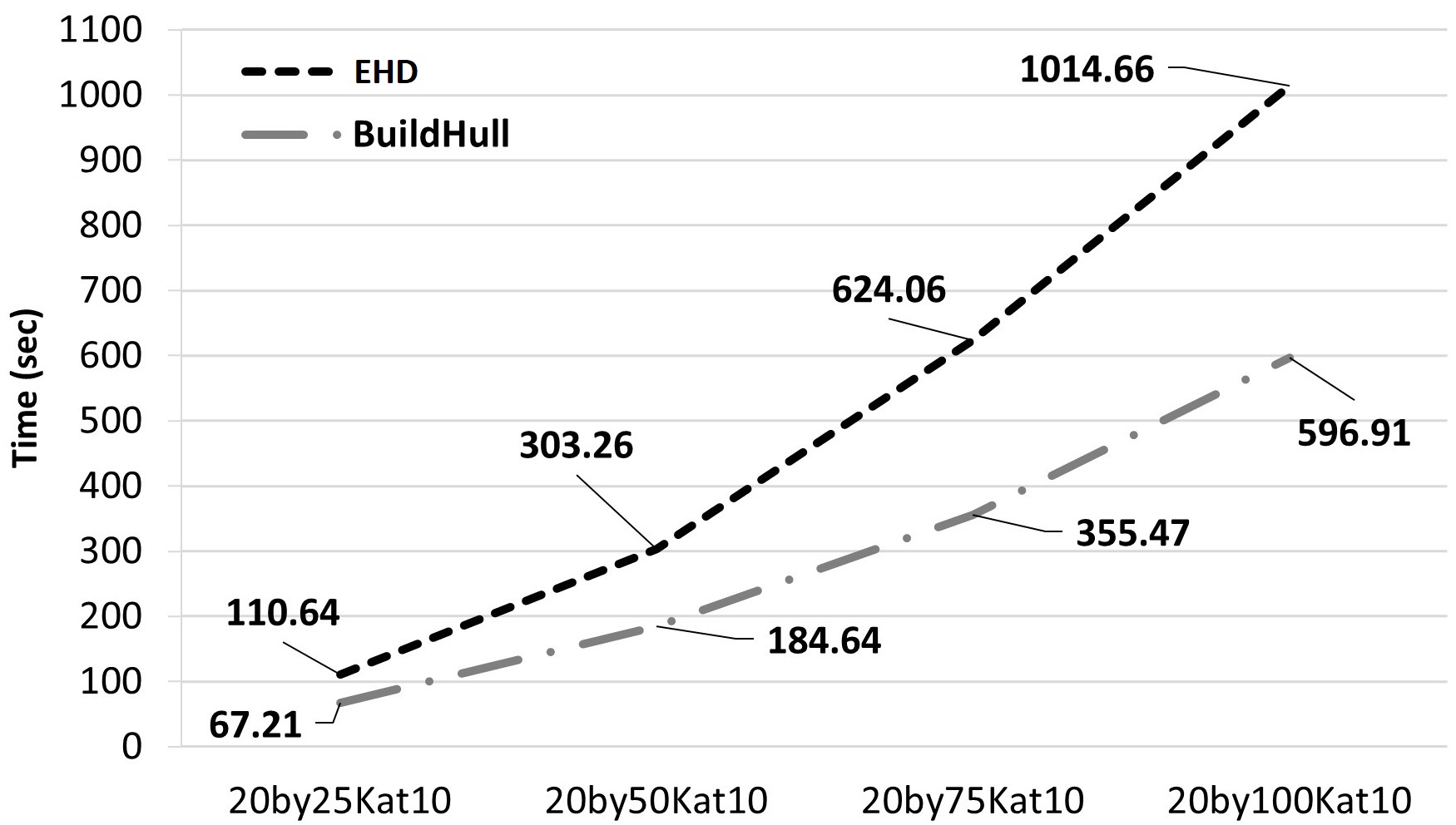}  
  \caption{$d=10\%$, $m=20$}
  %\label{fig:sub-second}
\end{subfigure}
\caption{Impact of cardinality on running time for the data sets in the {\tt MassiveScaleDEAdata} suite when the density is $d=10\%$.}
\label{fig:CardinalityD10}
\end{figure}

\noindent\underbar{Cardinality: $n$}. The four charts in Figure \ref{fig:CardinalityD10}  illustrate how  cardinality impacts execution time for the data sets with density $d=10\%$ --  a value for $d$ typical of what was observed throughout.
Similar to the effect of density, the running time performance gaps based on cardinality start small for the lowest cardinality, but widens as its value increases for all four dimensions.
For the experiments portrayed in Figure  \ref{fig:CardinalityD10}, the speedups of  \ttBH\  over \ttKZCT\ when $n=100$K range  from 1.3 to 2.3.
This performance gap does not grow uniformly by dimension $m$:  the ratio goes from 2.3 when $m=10$ to 1.6 for $m=15$.
We will see that data sets with  dimension $m=10$ will also be of interest when we study the effects of dimension next. 
%This seems to contradict what is expected and verified  about the role of dimension in DEA procedures, namely that as dimension increases, the time to process DEA data file takes longer. 
%%%12/30/21 This effect (the decrease in time when BOTH BH & KZCT go from m=10 to m=15 -- expecting an increase) needs to be explored/explained.
%%% GK Conjecture: the effect must be the preprocessors: DimSort, pre-scoring, & re-ordering based on pre-scoring.
%%% What happens in Dula11? GK generated the two graphs in Fig 2 
%%% What happens in Dula08?
%%% Can we discuss this why m=10 plays this role here and elsewhere? 
%%% Is there a "natural" disadvantage for KZCT in m=10?
%%% GK Conjecture: The pre-scoring affects the two procedures differently
%%% This is less consistent than what was observed when looking at effect of density.
The plots exhibit a familiar, although slight, quadratic behavior previously reported (see, e.g. \cite{Dula11}, \cite{Dula08}) for both procedures.

%%% The impact of dimension on running time for the three densities, d=1\%, d=10\%, d=25\%, with the four cardinalities for the two procedures are depicted in Figure \ref{fig:DimensionImpact}.

% In contrast, in Dula (2011) the execution time increases as the dimension increases.
% The current implementation is different in two ways, bigger initial partial hull and specific order for processing the DMUs.

\begin{figure}[H]
\begin{subfigure}{.5\textwidth}
  \centering
  \includegraphics[width=1\linewidth, height=0.75\textwidth]{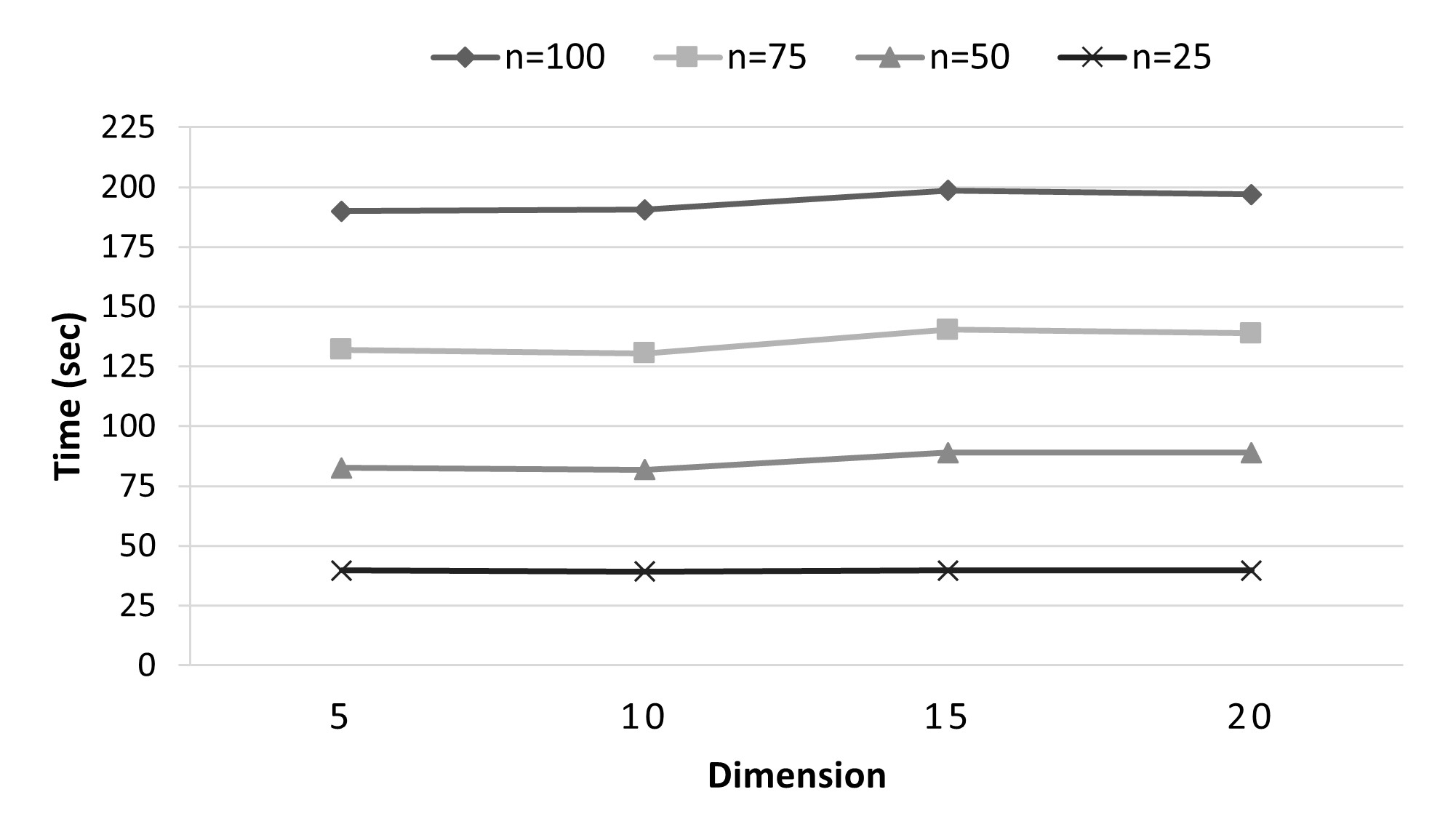}  
  \caption{\ttBH, $d=1\%$}
  \label{fig:BHfirst}
\end{subfigure}
\begin{subfigure}{.5\textwidth}
  \centering
  \includegraphics[width=1\linewidth, height=0.75\textwidth]{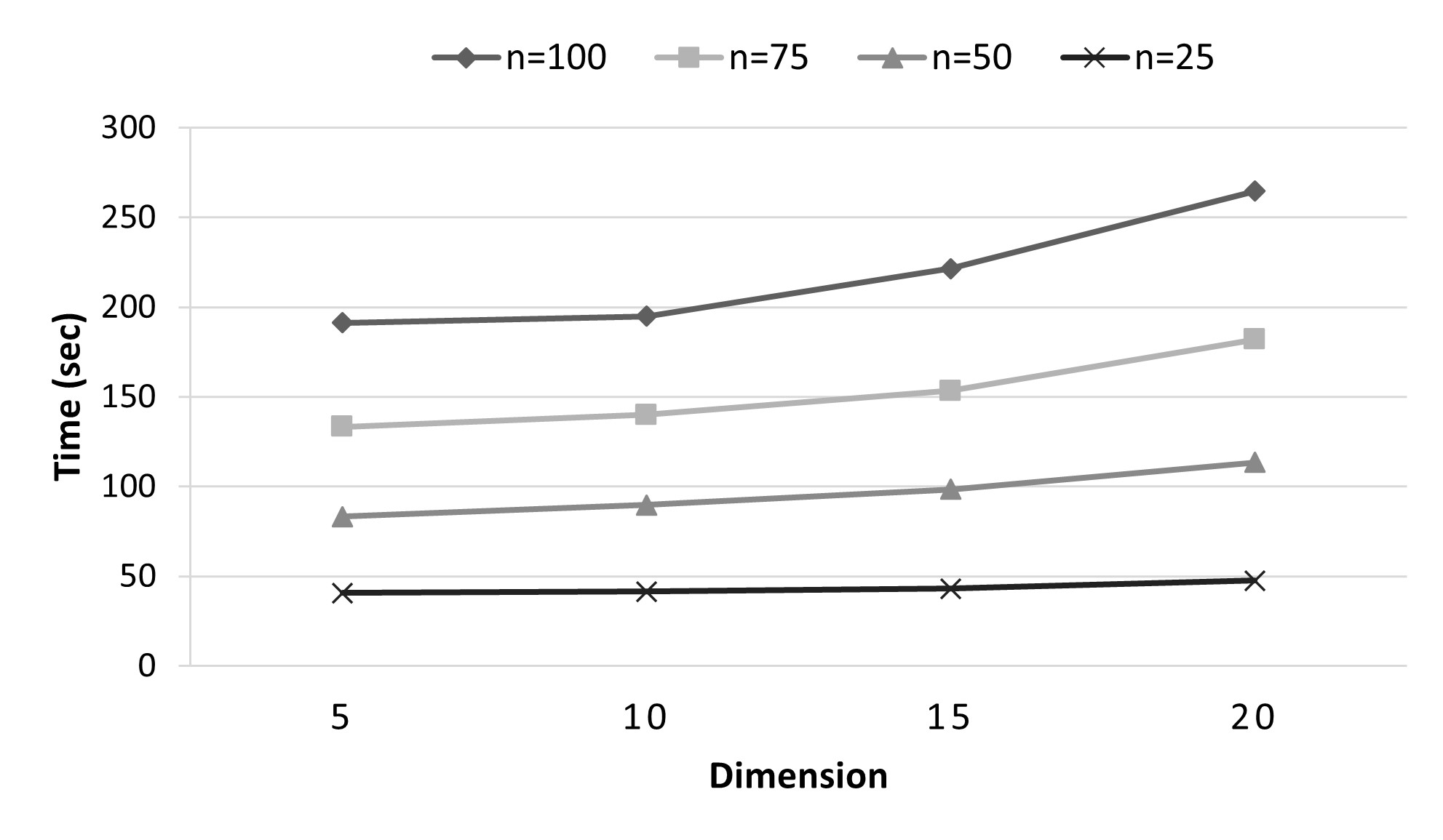}  
  \caption{\ttKZCT, $d=1\%$}
  \label{fig:KZCTfirst}
\end{subfigure}
\begin{subfigure}{.5\textwidth}
  \centering
  \includegraphics[width=1\linewidth, height=0.75\textwidth]{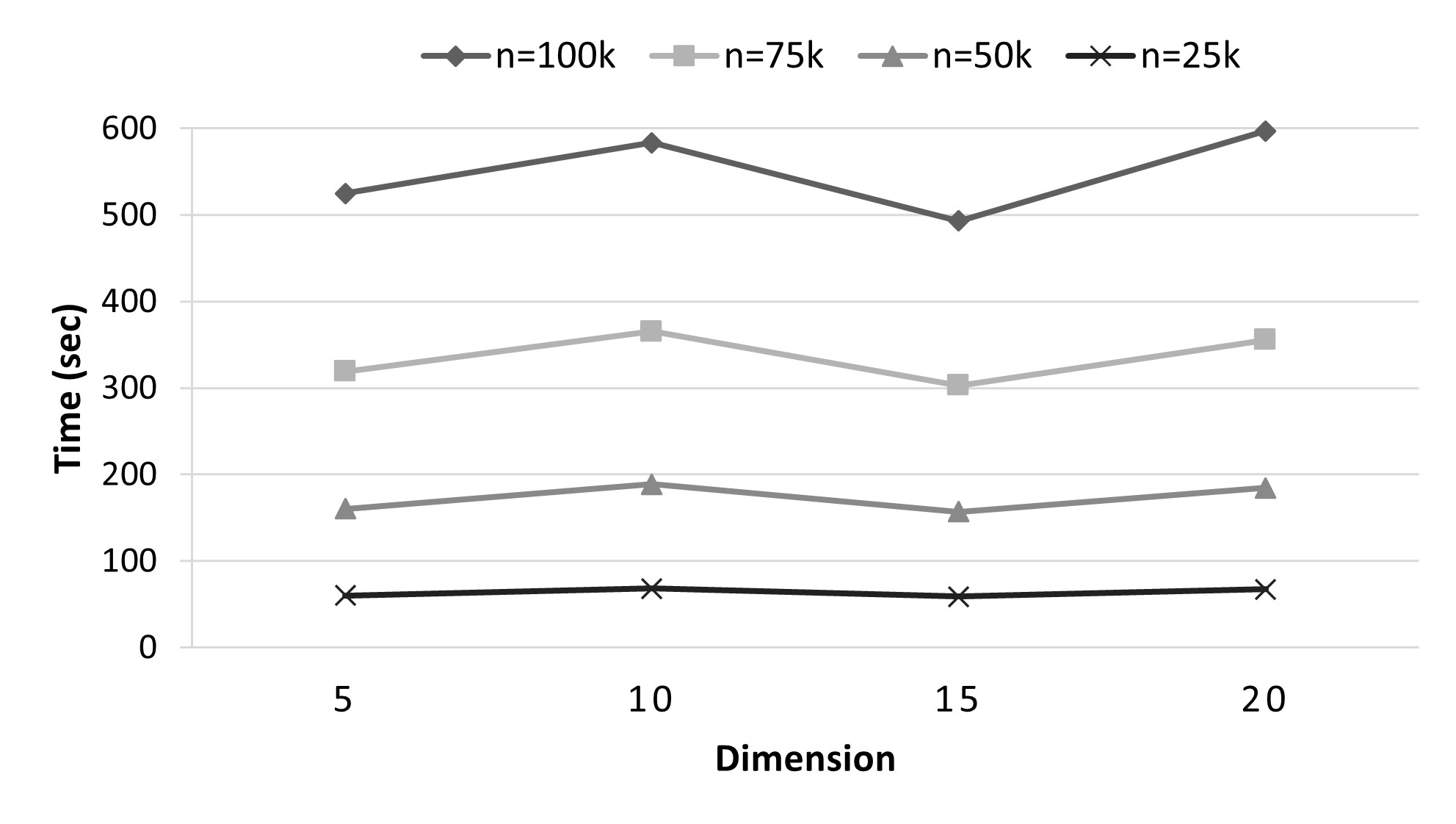}  
  \caption{\ttBH, $d=10\%$}
  \label{fig:BHsecond}
\end{subfigure}
\begin{subfigure}{.5\textwidth}
  \centering
  \includegraphics[width=1\linewidth, height=0.75\textwidth]{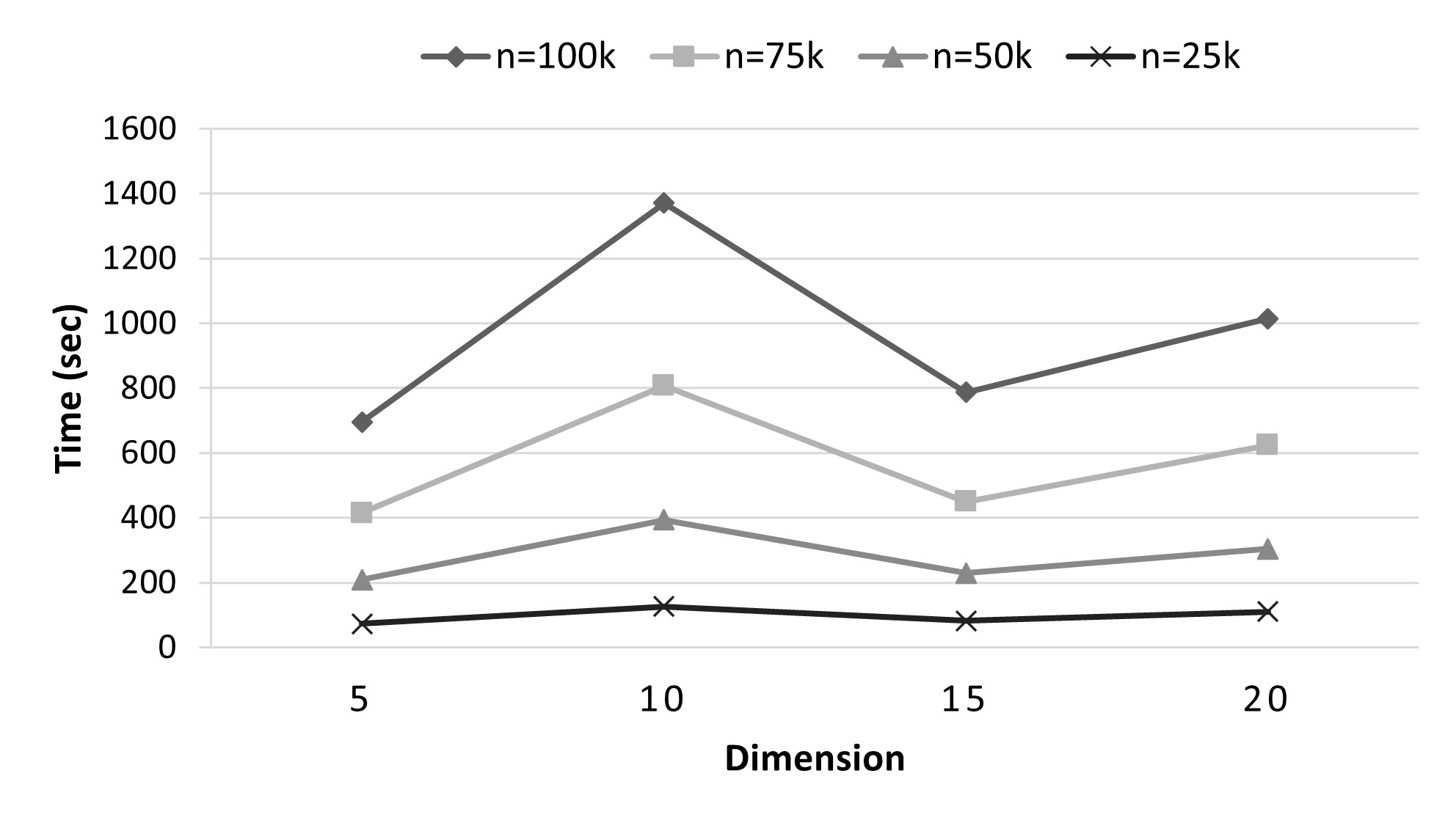}  
  \caption{\ttKZCT, $d=10\%$}
  \label{fig:KZCTsecond}
\end{subfigure}
\begin{subfigure}{.5\textwidth}
  \centering
  \includegraphics[width=1\linewidth, height=0.75\textwidth]{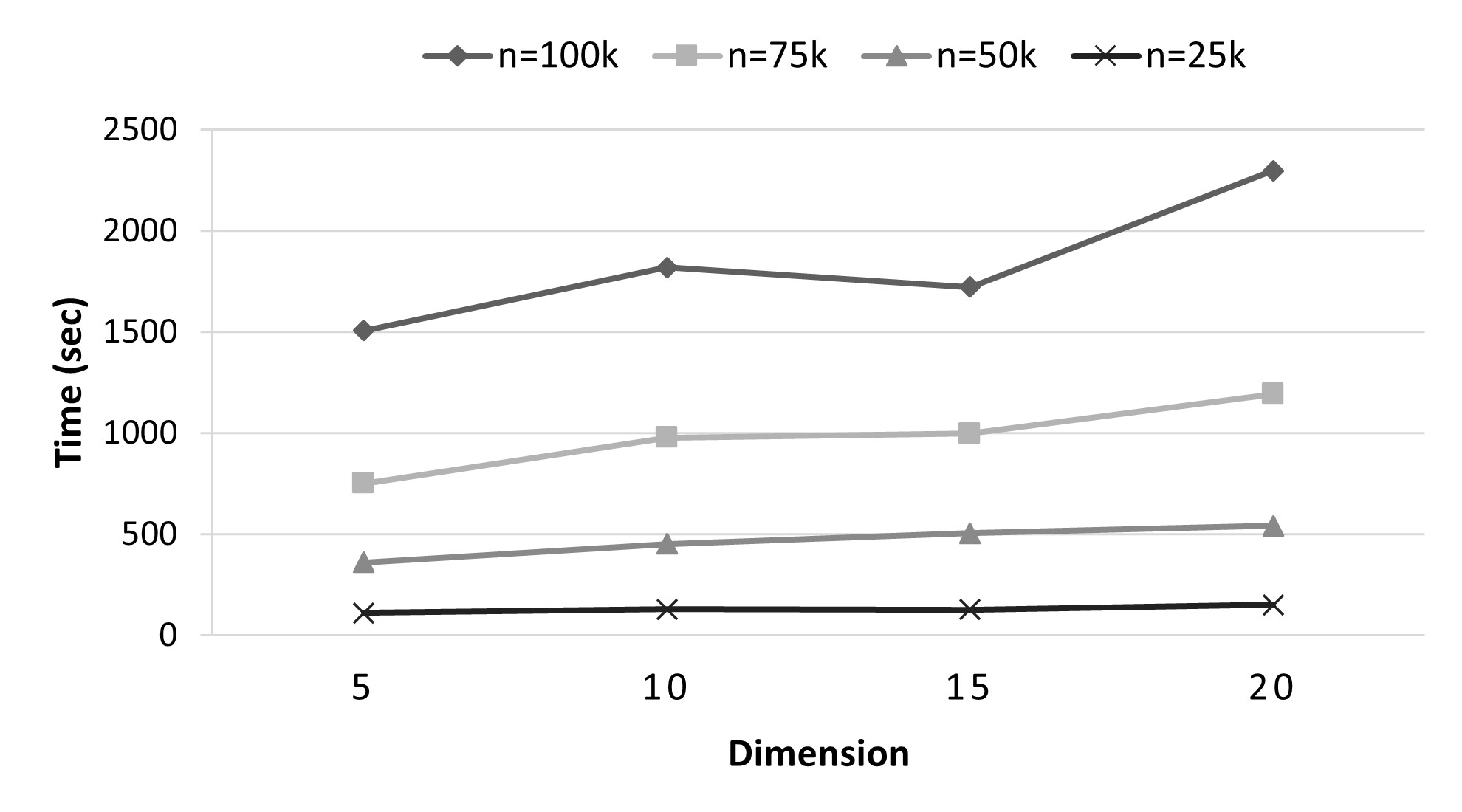}  
  \caption{\ttBH, $d=25\%$}
  \label{fig:BHthird}
\end{subfigure}
\begin{subfigure}{.5\textwidth}
  \centering
  \includegraphics[width=1\linewidth, height=0.75\textwidth]{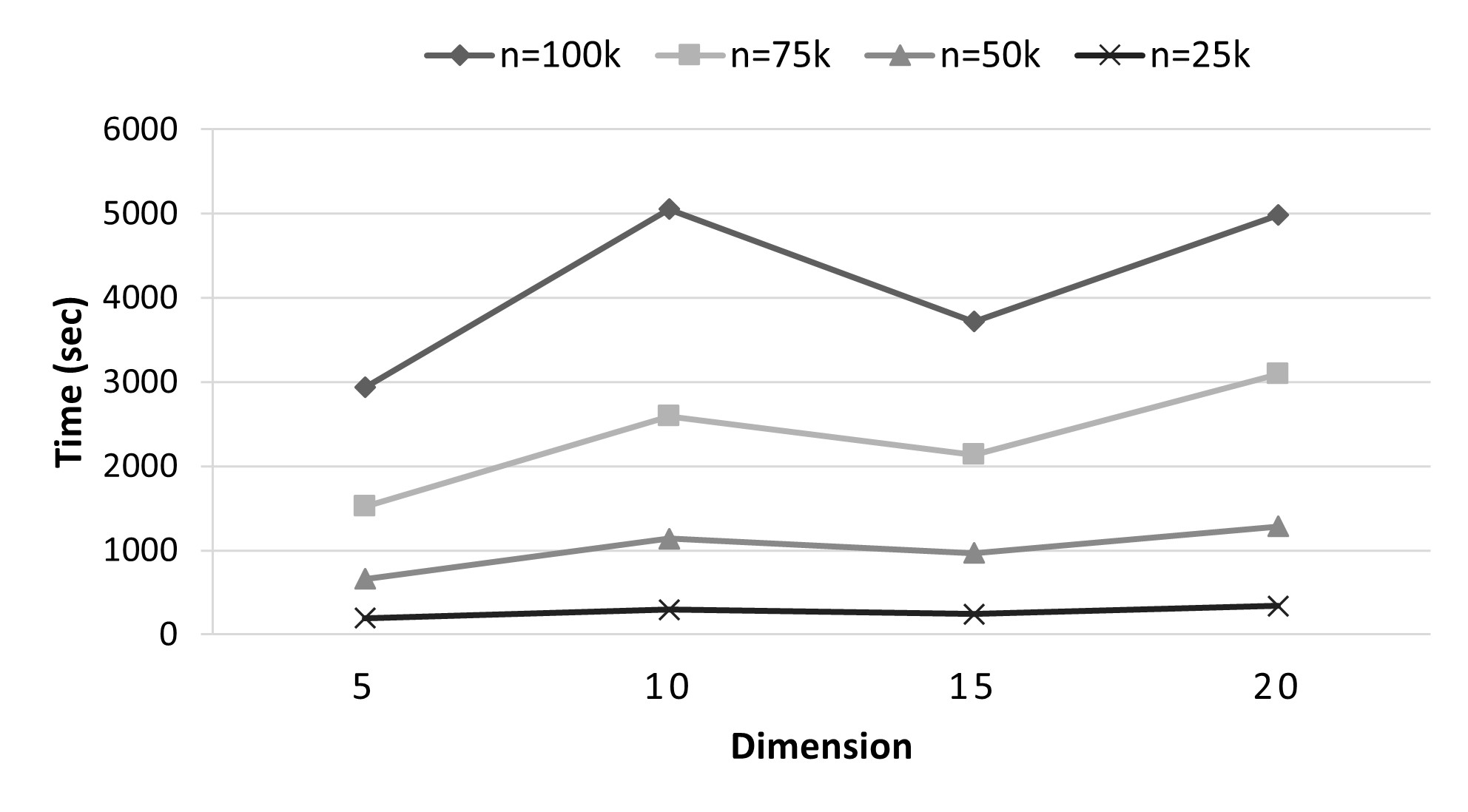}  
  \caption{\ttKZCT, $d=25\%$}
  \label{fig:KZCTthird}
\end{subfigure}
\caption{Impact of dimension on running time for data sets in the {\tt MassiveScaleDEAdata} suite.}
\label{fig:DimensionImpact}
\end{figure}

\par
\noindent\underbar{Dimension: $m$}.  The six panels in Figure \ref{fig:DimensionImpact} report on the relation between dimension $m$ and running time for the three densities $d=1\%, d=10\%, d=25\%$ in the \break  {\tt MassiveScaleDEAdata} suite for both procedures.
Each panel superimposes the measurements for the four cardinalities.
An unexpected reduction in running time from $m=10$ to $m=15$ for the four cardinalities is clearly visible in the panels (4c) \& (4d) and (4e) \& (4f) of Figure \ref{fig:DimensionImpact}.
This effect becomes more obvious as cardinality increases and more evident for \ttKZCT.
A similar effect was observed in \cite{Dula11} on the same {\tt MassiveScaleDEAdata} suite.
Results from \cite{Dula11} include instances of non-monotonicity, or a milder non-uniform rate of increase, in running time as dimension increases.
This occurs for data sets in the suite with $m=10$, e.g. {\tt **by25000at25} and {\tt **by100000at10} (see appendix in \cite{Dula11}).
The fact that the reduction in running time occurs for similar values of $m$ \& $d$ in both procedures in the current implementation is evidence that the phenomenon is connected to an intrinsic characteristic of the data sets involved.
In particular, for these datasets, the number and size of LPs solved in Step 4 of \ttKZCT\ are greater for $m=10$ than the corresponding ones with $m=15$.
For instance, from Table \ref{table:KZCT19Steps100K}, we see that, for data set {\tt 10by100000at10}, Step 4 of \ttKZCT\ solves 17,125 LPs each with 17,136 variables for an execution time of 1370.6 secs., whereas for the higher dimensional data set, {\tt 15by100000at10}, 10,097 LPs are solved each with 10,113 variables for an execution time of 787.18 secs.
For \ttBH, for these same two datasets, the average LP size is higher (1387.13) for $m=10$ than the one where $m=15$ (804.52), resulting in execution times of 583.48 secs. and 492.87 secs., respectively (see Table \ref{table:BestWorstBH100K}).

\par 
The non-intuitive effect where the increased dimension causes a reduction in execution time in some data sets appears in both the current  implementation of \ttBH\ and the one in 2011 \cite{Dula11}.
This suggests that the cause is related to characteristics of the data and how they affect the geometry of the partial hulls generated in the procedures. 
However, this also suggests that preprocessors are also playing a role.
In the case of \ttKZCT, the productivity of preprocessors is low for the eight datasets with $m=10$ dimensions, i.e., {\tt 10by$-n-$at10} and {\tt 10by$-n-$at25}; for $n=25K, 50K, 75K, 100K$, see Tables \ref{table:KZCT19Steps25K}-\ref{table:KZCT19Steps100K}.
Recall that we defined productivity for \ttKZCT\ as the number of interior data points from $\calA$ contained in the partial hull $\vrs(\calA^S)$ and identified in Steps 2 \& 3.
Productivity is a direct measurement of the effectiveness in Procedure \ttKZCT\ of the two preprocessors: dimension sorting and the pre-scoring heuristic.
Also, productivity determines the number and size of LPs solved in Step 4 of \ttKZCT: the higher the productivity of $\vrs(\calA^S)$, the closer the size of the LPs will be to the lower bound of $\vert\calF\vert$.
For instance, the productivity of {\tt 10by100000at10} dataset is almost 100\% while for {\tt 15by100000at10} is about 92\%.
The low productivity of the preprocessors for the aforementioned eight data sets explains the peculiar performance portrayed in panels (\ref{fig:KZCTsecond}) and (\ref{fig:KZCTthird}) of Figure \ref{fig:DimensionImpact}, as more LPs with larger size will be solved by \ttKZCT\ in these cases.

\par
Notice how \ttBH\ also feels the impact, though much more attenuated, of less effective preprocessing of these data sets as depicted in panels (\ref{fig:BHsecond}) and (\ref{fig:BHthird}) of Figure \ref{fig:DimensionImpact}.
These out-of-the-ordinary data sets prove to be useful, as they provide information about how the preprocessors impact the performance of the procedures.

\par
It is clear from the above analysis that preprocessors have a substantial impact on \ttKZCT\ and, when applied in an equitable manner to \ttBH, also renders the latter faster.
Ordering the rest of the points in ascending pre-score value facilitates to slow up the rate of growth of the LPs size in \ttBH.
The employed preprocessors work best under low density and high dimensionality.
Also, Step 4 in Procedure \ttKZCT\ requires considerable computational effort especially in instances where the cardinality of the frame is much larger than the paramenter $p$. 
This is verified by the experimentaion on {\tt MassiveScaleDEAdata} suite when the cardinality of the data set is large and its density is high.

\section{Conclusions}
Direct sequential (non-parallel) implementations of Procedures \ttKZCT\ from \cite{KhezrimotlaghZhuCookToloo19} and \ttBH\ in \cite{Dula11} were compared to obtain comprehensive empirical insights about their performances.
The comparisons are based on execution times but also using measurements about the number and size of the LPs solved to explain interesting performance characteristics.
The implementations establish a common ground for comparisons by applying the same preprocessors or enhancements when they fit.
The two procedures were tested on the publicly available and widely used {\tt MassiveScaleDEAdata} suite of 48 structured, synthetic, data sets.
The empirical results of these sequential implementations reveal that \ttBH\ is always faster than \ttKZCT, as speculated in \cite{Khezrimotlagh21}.
The computational time advantage of \ttBH\ relative to \ttKZCT\ is small for data sets with low density and increases as cardinality and density increase reaching speed ups of almost three times faster than \ttKZCT\  for the largest data sets.
Along with execution times, the total number of LPs solved was tracked and their size is calculated.
These metrics were used to explain interesting results such as how running time is reduced as dimension increases with both procedures.

\par 
A logical next step is a comparison of parallel implementations of \ttBH\ and \ttKZCT. 
Procedure \ttKZCT\ is directly parallelizable and this has been done in \cite{KhezrimotlaghZhuCookToloo19}.
The parallelization of \ttBH, although less immediate, will have to focus on more granular aspects of the algorithm such as classifying test points as external or internal to a partial hull.
The current work lays the foundation for a more theoretical study of the computational complexity of DEA procedures, something which is lacking in \ttBH\ and \ttKZCT\ and also in all the other important algorithmic contributions in DEA.
Another direction for future research is the hybridization of \ttBH\ with \ttKZCT, by applying the former separately to Steps 2, 3 \& 4 of the latter.
In addition, the insights obtained here about the role of enhancements and preprocessors  can be used to improve or develop variations, or new preprocessors outright.

\par 
Finally, contributions to computational DEA apply directly to the {\it frame} problem in computational geometry \cite{DulaLopez12}. 
Computational geometry offers a multitude of opportunities for applications of results from computational DEA in the form of algorithmic adaptations and comparisons along with applications on specialized problems outside of DEA.  

\appendix
\setcounter{secnumdepth}{0}
\section{Appendix}

\begin{table}[H]
  \centering
  \resizebox{\textwidth}{!}{
    \begin{tabular}{lcccccc}
    \hline
    Data Set & 
    \begin{tabular}[c]{@{}c@{}} LP Size \\ Step 2 \end{tabular} &
   \begin{tabular}[c]{@{}c@{}} LP Size \\ Step 3 \end{tabular} &
   \begin{tabular}[c]{@{}c@{}} LP Size \\ Step 4 \end{tabular} &
    \begin{tabular}[c]{@{}c@{}}\# of LPs \\ Step 4 \end{tabular} &
    \begin{tabular}[c]{@{}c@{}}Total \# \\ of LPs \end{tabular} &
    \begin{tabular}[c]{@{}c@{}}Total \\ Time \end{tabular} \\
    \hline
    05by25000at01 & 159   & 42    & 926   & 920   & 25915  & 40.70 \\
    05by25000at10 & 159   & 151   & 3304  & 3298  & 28293  & 73.03 \\
    05by25000at25 & 159   & 159   & 7247  & 7241  & 32236  & 192.43 \\
    10by25000at01 & 159   & 81    & 262   & 251   & 25241  & 41.36 \\
    10by25000at10 & 159   & 140   & 4602  & 4591  & 29581  & 126.34 \\
    10by25000at25 & 159   & 155   & 8119  & 8108  & 33098  & 296.17 \\
    15by25000at01 & 159   & 96    & 503   & 487   & 25472  & 43.06 \\
    15by25000at10 & 159   & 159   & 2608  & 2592  & 27577  & 81.95 \\
    15by25000at25 & 159   & 159   & 6382  & 6366  & 31351  & 244.25 \\
    20by25000at01 & 159   & 138   & 262   & 241   & 25221  & 47.55 \\
    20by25000at10 & 159   & 159   & 3045  & 3024  & 28004  & 110.64 \\
    20by25000at25 & 159   & 159   & 6683  & 6662  & 31642  & 337.29 \\
    \hline
    \end{tabular}
    }
   \caption{Results of \ttKZCT\ for cardinality $n=25$K in {\tt MassiveScaleDEAdata} suite.}
   \label{table:KZCT19Steps25K}
\end{table}

\begin{table}[H]
  \centering
  \resizebox{\textwidth}{!}{
\begin{tabular}{lcccccc}
    \hline
    Data Set & 
    \begin{tabular}[c]{@{}c@{}} LP Size \\ Step 2 \end{tabular} &
   \begin{tabular}[c]{@{}c@{}} LP Size \\ Step 3 \end{tabular} &
   \begin{tabular}[c]{@{}c@{}} LP Size \\ Step 4 \end{tabular} &
    \begin{tabular}[c]{@{}c@{}}\# of LPs \\ Step 4 \end{tabular} &
    \begin{tabular}[c]{@{}c@{}}Total \# \\ of LPs \end{tabular} &
    \begin{tabular}[c]{@{}c@{}}Total \\ Time \end{tabular} \\
    \hline
    05by50000at01 & 224   & 71    & 1596  & 1590  & 51585  & 83.14 \\
    05by50000at10 & 224   & 222   & 6531  & 6525  & 56520  & 208.53 \\
    05by50000at25 & 224   & 224   & 14643 & 14637 & 64632  & 660.96 \\
    10by50000at01 & 224   & 155   & 502   & 491   & 50481  & 89.68 \\
    10by50000at10 & 224   & 202   & 8887  & 8876  & 58866  & 393.12 \\
    10by50000at25 & 224   & 220   & 15869 & 15858 & 65848  & 1140.77 \\
    15by50000at01 & 224   & 167   & 619   & 603   & 50588  & 98.41 \\
    15by50000at10 & 224   & 224   & 5093  & 5077  & 55062  & 228.86 \\
    15by50000at25 & 224   & 224   & 12638 & 12622 & 62607  & 967.19 \\
    20by50000at01 & 224   & 210   & 506   & 485   & 50465  & 113.48 \\
    20by50000at10 & 224   & 224   & 5584  & 5563  & 55543  & 303.26 \\
    20by50000at25 & 224   & 224   & 13141 & 13120 & 63100  & 1283.35 \\
    \hline
    \end{tabular}
    }
   \caption{Results of \ttKZCT\ for cardinality $n=50$K in {\tt MassiveScaleDEAdata} suite.}
\label{table:KZCT19Steps50K}
\end{table}

\begin{table}[H]
  \centering
  \resizebox{\textwidth}{!}{
\begin{tabular}{lcccccc}
    \hline
    Data Set & 
    \begin{tabular}[c]{@{}c@{}} LP Size \\ Step 2 \end{tabular} &
   \begin{tabular}[c]{@{}c@{}} LP Size \\ Step 3 \end{tabular} &
   \begin{tabular}[c]{@{}c@{}} LP Size \\ Step 4 \end{tabular} &
    \begin{tabular}[c]{@{}c@{}}\# of LPs \\ Step 4 \end{tabular} &
    \begin{tabular}[c]{@{}c@{}}Total \# \\ of LPs \end{tabular} &
    \begin{tabular}[c]{@{}c@{}}Total \\ Time \end{tabular} \\
    \hline
    05by75000at01 & 274   & 84    & 2561  & 2555  & 77550  & 133.20 \\
    05by75000at10 & 274   & 269   & 10078 & 10072 & 85067  & 414.42 \\
    05by75000at25 & 274   & 274   & 22175 & 22169 & 97164  & 1522.59 \\
    10by75000at01 & 274   & 197   & 754   & 743   & 75733  & 139.96 \\
    10by75000at10 & 274   & 247   & 13001 & 12990 & 87980  & 806.96 \\
    10by75000at25 & 274   & 272   & 23988 & 23977 & 98967  & 2592.52 \\
    15by75000at01 & 274   & 202   & 827   & 811   & 75796  & 153.39 \\
    15by75000at10 & 274   & 274   & 7600  & 7584  & 82569  & 449.36 \\
    15by75000at25 & 274   & 274   & 18891 & 18875 & 93860  & 2134.45 \\
    20by75000at01 & 274   & 263   & 752   & 731   & 75711  & 181.88 \\
    20by75000at10 & 274   & 274   & 8219  & 8198  & 83178  & 624.06 \\
    20by75000at25 & 274   & 274   & 19419 & 19398 & 94378  & 3092.35 \\
    \hline
    \end{tabular}
    }
   \caption{Results of \ttKZCT\ for cardinality $n=75$K in {\tt MassiveScaleDEAdata} suite.}
\label{table:KZCT19Steps75K}
\end{table}

\begin{table}[H]
  \centering
  \resizebox{\textwidth}{!}{
\begin{tabular}{lcccccc}
    \hline
    Data Set & 
    \begin{tabular}[c]{@{}c@{}} LP Size \\ Step 2 \end{tabular} &
   \begin{tabular}[c]{@{}c@{}} LP Size \\ Step 3 \end{tabular} &
   \begin{tabular}[c]{@{}c@{}} LP Size \\ Step 4 \end{tabular} &
    \begin{tabular}[c]{@{}c@{}}\# of LPs \\ Step 4 \end{tabular} &
    \begin{tabular}[c]{@{}c@{}}Total \# \\ of LPs \end{tabular} &
    \begin{tabular}[c]{@{}c@{}}Total \\ Time \end{tabular} \\
    \hline
    05by100000at01 & 317   & 88    & 3477  & 3471  & 103466   & 191.19 \\
    05by100000at10 & 317   & 315   & 13636 & 13630 & 113625   & 694.81 \\
    05by100000at25 & 317   & 317   & 29519 & 29513 & 129508   & 2937.29 \\
    10by100000at01 & 317   & 225   & 1001  & 990   & 100980   & 194.96 \\
    10by100000at10 & 317   & 296   & 17136 & 17125 & 117115   & 1370.60 \\
    10by100000at25 & 317   & 314   & 32178 & 32167 & 132157   & 5051.58 \\
    15by100000at01 & 317   & 253   & 1108  & 1092  & 101077   & 221.51 \\
    15by100000at10 & 317   & 317   & 10113 & 10097 & 110082   & 787.18 \\
    15by100000at25 & 317   & 317   & 25148 & 25132 & 125117   & 3716.97 \\
    20by100000at01 & 317   & 308   & 1005  & 984   & 100964   & 264.65 \\
    20by100000at10 & 317   & 317   & 10762 & 10741 & 110721   & 1014.66 \\
    20by100000at25 & 317   & 317   & 25830 & 25809 & 125789   & 4983.47 \\
    \hline
    \end{tabular}
    }
  \caption{Results of \ttKZCT\ for cardinality $n=100$K in {\tt MassiveScaleDEAdata} suite.}
\label{table:KZCT19Steps100K}
\end{table}

\begin{table}[H]
  \centering
  %\resizebox{\textwidth}{!}{
    \begin{tabular}{lcccc}
    \hline
    Data Set & 
    \begin{tabular}[c]{@{}c@{}}Total \# \\ of LPs \end{tabular}&
    \begin{tabular}[c]{@{}c@{}}Average \\ LP Size \end{tabular} &    
    \begin{tabular}[c]{@{}c@{}}Total \\ Time 
    \end{tabular}&
    \\
    \hline
    05by25000at01  & 24995 & 45.74 & 39.74\\
    05by25000at10  & 24995 & 424.12 & 60.27 \\
    05by25000at25  & 24995 & 1331.31 & 110.98 \\
    10by25000at01  & 24990 & 33.51 & 39.18\\
    10by25000at10  & 24990 & 429.01 & 68.18 \\
    10by25000at25  & 24990 & 1209.34 & 129.13 \\
    15by25000at01  & 24985 & 66.22 & 39.71\\
    15by25000at10  & 24985 & 221.18 & 59.08 \\
    15by25000at25  & 24985 & 908.63 & 127.19 \\
    20by25000at01  & 24980 & 48.91 & 39.78 \\
    20by25000at10  & 24980 & 236.70  & 67.21\\
    20by25000at25  & 24980 & 905.09 & 151.02 \\
    \hline
    \end{tabular}
   % }
    \caption{Results of \ttBH\ for cardinality $n=25$K in {\tt MassiveScaleDEAdata} suite.}
\label{table:BestWorstBH25K}
\end{table}

\begin{table}[H]
  \centering
  %\resizebox{\textwidth}{!}{
        \begin{tabular}{lcccc}
    \hline
    Data Set & 
    \begin{tabular}[c]{@{}c@{}}Total \# \\ of LPs \end{tabular}&
    \begin{tabular}[c]{@{}c@{}}Average \\ LP Size  \end{tabular} &   
    \begin{tabular}[c]{@{}c@{}}Total \\ Time\end{tabular}&
 \\
    \hline
    05by50000at01  & 49995 & 74.31 & 82.65 \\
    05by50000at10  & 49995 & 818.03  & 160.48\\
    05by50000at25  & 49995 & 2638.61 & 359.40 \\
    10by50000at01  & 49990 & 45.96 & 81.79\\
    10by50000at10  & 49990 & 734.84  & 188.57\\
    10by50000at25  & 49990 & 2380.06 & 452.18 \\
    15by50000at01  & 49985 & 74.92  & 88.96\\
    15by50000at10 & 49985 & 413.03 & 156.90  \\
    15by50000at25 & 49985 & 1793.22  & 503.64 \\
    20by50000at01  & 49980 & 52.65 & 89.02  \\
    20by50000at10  & 49980 & 412.23 & 184.64 \\
    20by50000at25  & 49980 & 1753.61 & 541.08\\
    \hline
    \end{tabular}
    %}
    \caption{Results of \ttBH\ for cardinality $n=50$K in {\tt MassiveScaleDEAdata} suite.}
\label{table:BestWorstBH50K}
\end{table}

\begin{table}[H]
  \centering
 % \resizebox{\textwidth}{!}{
    \begin{tabular}{lcccc}
    \hline
    Data Set & 
    \begin{tabular}[c]{@{}c@{}}Total \# \\ of LPs \end{tabular}&
    \begin{tabular}[c]{@{}c@{}}Average \\ LP Size  \end{tabular} &   
    \begin{tabular}[c]{@{}c@{}}Total \\ Time\end{tabular}&
 \\
    \hline
    05by75000at01  & 74995 & 114.66 & 131.94 \\
    05by75000at10  & 74995 & 1241.59 & 318.91 \\
    05by75000at25  & 74995 & 3924.94 & 750.39 \\
    10by75000at01  & 74990 & 62.53 & 130.42 \\
    10by75000at10  & 74990 & 1073.09 & 365.38 \\
    10by75000at25  & 74990 & 3617.11  & 977.13\\
    15by75000at01  & 74985 & 82.95 & 140.36\\
    15by75000at10  & 74985 & 621.01 & 303.11  \\
    15by75000at25  & 74985 & 2710.10 & 997.74 \\
    20by75000at01 & 74980 & 60.54 & 138.79  \\
    20by75000at10  & 74980 & 629.17 & 355.47\\
    20by75000at25  & 74980 & 2598.93 & 1192.76\\
    \hline
    \end{tabular}
   % }
    \caption{Results of \ttBH\ for cardinality $n=75$K in {\tt MassiveScaleDEAdata} suite.}
\label{table:BestWorstBH75K}
\end{table}

\begin{table}[H]
  \centering
  %\resizebox{\textwidth}{!}{
   \begin{tabular}{lcccc}
    \hline
    Data Set & 
    \begin{tabular}[c]{@{}c@{}}Total \# \\ of LPs \end{tabular}&
    \begin{tabular}[c]{@{}c@{}}Average \\ LP Size  \end{tabular} &   
    \begin{tabular}[c]{@{}c@{}}Total \\ Time\end{tabular}&
 \\
    \hline
    05by100000at01  & 99995 & 149.74 & 189.87 \\
    05by100000at10  & 99995 & 1666.37 & 524.86 \\
    05by100000at25  & 99995 & 5228.89 & 1506.57 \\
    10by100000at01  & 99990 & 80.25 & 190.45 \\
    10by100000at10  & 99990 & 1387.13 & 583.48 \\
    10by100000at25  & 99990 & 4676.17 & 1817.77  \\
    15by100000at01  & 99985 & 89.91 & 198.66 \\
    15by100000at10  & 99985 & 804.52 & 492.87 \\
    15by100000at25  & 99985 & 3542.93 & 1721.24 \\
    20by100000at01  & 99980 & 64.10 & 196.85 \\
    20by100000at10  & 99980 & 821.44 & 596.91 \\
    20by100000at25  & 99980 & 3561.68 & 2295.94\\
    \hline
    \end{tabular}
    %}
    \caption{Results of \ttBH\ for cardinality $n=100$K in {\tt MassiveScaleDEAdata} suite.}
\label{table:BestWorstBH100K}
\end{table}

\bibliographystyle{plain} % outcomment this and next line in Case 1
\bibliography{References} % if more than one, comma separated
%%%%%%%%%%%%%%%%%%%%%%%%%%%%%%%%%%%%%%%%%%%%%
\end{document}